\documentclass{compositio}
\usepackage{latexsym}
\usepackage{amscd, amsfonts, eucal, mathrsfs, amsmath, amssymb, amsthm}
\usepackage[OT2,OT1]{fontenc}
\def\cyr{\fontencoding{OT2}\fontfamily{wncyr}\selectfont}
\input xy
\xyoption{all}

\newcommand{\Caldararu}{C\u ald\u araru}

\newcommand{\field}[1]{\mathbf #1}
\newcommand{\mf}[1]{\mathfrak #1}
\newcommand{\mc}[1]{\mathcal #1}
\newcommand{\ms}[1]{\mathscr #1}
\newcommand{\widebar}[1]{\overline{#1}}

\newcommand{\R}{\field R}

\newcommand{\C}{\field C}
\newcommand{\F}{\field F}
\newcommand{\Z}{\field Z}
\newcommand{\Q}{\field Q}

\newcommand{\simto}{\stackrel{\sim}{\to}}

\renewcommand{\phi}{\varphi}

\renewcommand{\O}{\ms O}

\newcommand{\shom}{\ms H\!om}

\newcommand{\rshom}{\mathbf{R}\shom}

\DeclareMathOperator{\sh}{Sh}

\newcommand{\send}{\ms E\!nd}

\newcommand{\spec}{\operatorname{Spec}}

\renewcommand{\P}{\field P}

\DeclareMathOperator{\Pic}{Pic}
\newcommand{\td}{\operatorname{Td}}
\DeclareMathOperator{\chern}{ch}

\DeclareMathOperator{\pr}{pr}

\DeclareMathOperator{\bl}{Bl}

\newcommand{\m}{\boldsymbol{\mu}}

\newcommand{\G}{\field G} 

\renewcommand{\H}{\operatorname{H}}

\newcommand{\GL}{\operatorname{GL}}
\newcommand{\PGL}{\operatorname{PGL}}
\DeclareMathOperator{\SL}{\operatorname{SL}}

\newcommand{\Sha}{\text{\cyr Sh}}

\DeclareMathOperator{\qcoh}{\operatorname{QCoh}}

\newcommand{\ch}{\operatorname{char}}
\DeclareMathOperator{\per}{per}
\DeclareMathOperator{\ind}{ind}

\DeclareMathOperator{\D}{\operatorname{\bf D}}

\DeclareMathOperator*{\tensor}{\otimes}
\DeclareMathOperator*{\ltensor}{\stackrel{\field L}{\otimes}}

\DeclareMathOperator{\rk}{\operatorname{rk}}

\newcommand{\surj}{\twoheadrightarrow}
\newcommand{\inj}{\hookrightarrow}

\DeclareMathOperator{\aut}{\operatorname{Aut}}
\DeclareMathOperator{\isom}{\operatorname{Isom}}
\DeclareMathOperator{\M}{\operatorname{M}}
\DeclareMathOperator{\Br}{\operatorname{Br}}

\newcommand{\Tw}{\mathbf{Tw}}
\DeclareMathOperator{\mTw}{Tw}
\newcommand{\Sh}{\mathbf{Sh}}
\DeclareMathOperator{\mSh}{Sh}
\DeclareMathOperator{\B}{\operatorname{B\!}}
\DeclareMathOperator{\Mod}{\operatorname{\ms M\!od}}

\newtheorem{lem}{Lemma}[subsubsection]
\newtheorem{LEM}{Lemma}[section]
\newtheorem{thm}[lem]{Theorem}
\newtheorem{THM}[LEM]{Theorem}
\newtheorem*{theorem}{Theorem}

\newtheorem{prop}[lem]{Proposition}
\newtheorem{PROP}[LEM]{Proposition}

\newtheorem{cor}[lem]{Corollary}
\newtheorem{COR}[LEM]{Corollary}

\newtheorem*{conj}{Conjecture}
\newtheorem{ques}[lem]{Question} 

\theoremstyle{definition}
\newtheorem{defn}[lem]{Definition}
\newtheorem{DEFN}[LEM]{Definition}

\newtheorem{example}[lem]{Example}

\newtheorem{para}[lem]{\indent}
\newtheorem{situation}[lem]{Situation}

\theoremstyle{remark}
\newtheorem{remark}[lem]{Remark}
\newtheorem{REMARK}[LEM]{Remark}

\newtheorem{notn}[lem]{Notation}

\begin{document}

\title{Twisted sheaves and the period-index problem}

\author{Max Lieblich}
\email{lieblich@math.princeton.edu}
\address{Department of Mathematics, Fine Hall, Washington Road,
  Princeton NJ 08544-1000}
\thanks{The author was partially supported by a Clay Liftoff
  Fellowship and an NSF Postdoctoral Fellowship during the preparation
of this paper.}

\bibliographystyle{amsalpha}

\begin{abstract} We use twisted sheaves and their moduli spaces to
  study the Brauer group of a scheme.  In particular, we (1) show how
  twisted methods can be efficiently used to re-prove the basic facts
  about the Brauer group and cohomological Brauer group (including
  Gabber's theorem that they coincide for a separated union of two
  affine schemes), (2) give a new proof of de Jong's period-index
  theorem for surfaces over algebraically closed fields, and (3) prove
  an analogous result for surfaces over finite fields.  We also
  include a reduction of all period-index problems for Brauer groups
  of function fields over algebraically closed fields to
  characteristic $0$, which (among other things) extends de Jong's
  result to include classes of period divisible by the characteristic
  of the base field.  Finally, we use the theory developed here to
  give counterexamples to a standard type of local-to-global conjecture
  for geometrically rational varieties over the function field of the
  projective plane.
\end{abstract}

\maketitle

\tableofcontents

\section{Introduction}

Let $K$ be a field.  Given a Brauer class $\alpha\in\Br(K)$, there are
two natural discrete invariants: the \emph{period\/} of $\alpha$,
$\per(\alpha)$, which is the order of $\alpha$ as an element of
$\Br(K)$, and the \emph{index\/} of $\alpha$, $\ind(\alpha)$, whose
square is the $K$-vector space dimension of a central division algebra
with Brauer class $\alpha$.  It is elementary to prove that
$\per(\alpha)|\ind(\alpha)$ and that they have the same prime factors,
so that $\ind(\alpha)|\per(\alpha)^{\ell}$ for some $\ell$ (see
\ref{L:basic relation} below).  The properties of the exponent $\ell$
have been studied by numerous authors
(\cite{artinperind},\cite{colliot-thelene},\cite{dejong-per-ind},\cite{guletskii},\cite{kresch},\cite{nakayama},\cite{saltman},\cite{vandenbergh},\ldots),
and the work to date has all supported the following folk-conjecture.

\begin{conj} If $K$ is a $C_d$ field and $\alpha\in\Br(K)$, then
  $\ind(\alpha)|\per(\alpha)^{d-1}$.
\end{conj}
 
The main goal of this paper is to prove the following theorem.

\begin{theorem} Suppose $K$ is a finitely generated field of
  characteristic $p>0$ with transcendence degree $2$ over the prime
  field.  If $\alpha\in\Br(K)$ has order prime to $\ch(K)$, then
  $\ind(\alpha)|\per(\alpha)^3$.  If $\alpha$ is unramified, then
  $\ind(\alpha)=\per(\alpha)$.
\end{theorem}

The reader will note that the exponent we obtain in the ramified case
should be $2$ according to the conjecture (although the unramified
result is somewhat better than predicted).  This gap is close to being
filled, at the expense of adding another layer of stack-theoretic
complexity to the methods employed here.

The class $\alpha$ is unramified precisely when it comes from a class
in the Brauer group of a smooth projective model $X$ for the function
field $K$.  As we will describe below, techniques of Saltman
\cite{saltman} permit the deduction of the ramified case from the
unramified one.  (It is precisely here where one is forced to increase
the exponent by $2$ rather than by $1$.)  The unramified case is
proven by finding a point in a certain moduli space parametrizing
solutions to the equation $\per=\ind$.  We must of course first show
that this moduli space is non-empty, which is precisely de Jong's
period-index theorem \cite{dejong-per-ind}.  We offer a new proof of
this result, as well as a reduction of all period-index problems for
function fields over algebraically closed fields to characteristic
$0$, thus giving an incremental improvement of de Jong's theorem: we
will not need the hypothesis that the period of $\alpha$ is prime to
the characteristic of the base field (until we start working over
finite base fields).

The moduli space we will use to prove these theorems is the space of
\emph{stable twisted sheaves\/}.  As this paper is intended to
encourage the use of twisted methods in the study of the Brauer group,
we have included more than is strictly necessary for our purposes.  We
give a slightly new proof of Gabber's theorem that the Brauer group
and cohomological Brauer group of an affine scheme coincide.  (For an
introduction to this problem, as well as a proof of Gabber's theorem
of which ours is an adaptation, the reader is referred to
\cite{hoobler}.)  The methods used here are a precursor to de Jong's
proof of Gabber's more recent theorem: $\Br=\Br'$ for any
quasi-compact separated scheme admitting an ample invertible sheaf
\cite{dejong-gabber}.

The philosophy of twisted sheaves leads naturally to a question
(\ref{Conj:conj} below) about the geometry of the moduli space of
stable vector bundles with fixed determinant on a curve, whose answer
is closely related to the period-index problem (in any dimension).  By
relating these two questions, we are able to give numerous examples of
smooth projective geometrically connected geometrically rational
varieties with geometric Picard group $\Z$ over $\C(t_1,t_2)$ which
have points over the completion of $\C(t_1,t_2)$ at any (arbitrary)
valuation but which lack rational points.  These appear to be the
first examples of such varieties over $\C(t_1,t_2)$ which violate the
Hasse principle.  (Hassett has informed me that there are also
unpublished examples of cubic surfaces over $\C(t_1,t_2)$ which
violate the Hasse principle and have vanishing Brauer-Manin
obstruction; these of course have higher geometric Picard number.)  It
seems that the method given here for producing such examples might be
useful in other contexts.

A brief summary of the contents of this paper: In section
\ref{S:background} we sketch the basic aspects of the period-index
problem and the theory of gerbes (apr\`es Giraud); the reader familiar
with this material can skip this section, referring back for
terminology as needed.  This is followed in section \ref{S:twisted
  sheaves and their moduli} by a rapid tour of the theory of twisted
sheaves, their use in proving the basic facts about the Brauer group
of a scheme and Gabber's theorem, and a summary of the theory of their
moduli on low-dimensional ambient varieties (fully developed in
\cite{twisted-moduli}).  Finally, in section \ref{S:period-index} we
apply the abstract theory to the period-index problem and prove the
main theorem.  We have included a very short appendix discussing a few
basic facts about (quasi-)coherent sheaves on stacks and elementary
transforms.

\begin{acknowledgements}
  This paper (along with \cite{twisted-moduli} and \cite{pgl-bundles})
  is based upon part of my Ph.D. thesis, written under the guidance of
  Aise Johan de Jong \cite{mythesis}.  I also received numerous
  helpful comments from Jean-Louis Colliot-Th\'el\`ene, Laurent
  Moret-Bailly, Martin Olsson, David Saltman, Jason Starr, Burt
  Totaro, and Olivier Wittenberg.  Finally, I owe a great debt of
  gratitude to the referee, who gave this paper a very careful reading
  and helped to drastically improve the exposition.  As usual, any
  remaining errors are entirely my own.
\end{acknowledgements}

\section{Background}\label{S:background}

\subsection{Algebra}

In this section we rapidly summarize the basic facts concerning the
period-index problem.

\subsubsection{Period and index}

The period-index problem has several manifestations.  One way of
describing the problem uses Galois cohomology.  Let $K$ be a field and
$F$ an \'etale sheaf on $\spec K$.  Let $\alpha\in\H^{i}(\spec K,F)$
be a Galois cohomology class with $i>0$.

\begin{defn}\label{D:per-ind defn} 
  The \emph{period\/} of $\alpha$, denoted $\per(\alpha)$, is the
  order of $\alpha$ in the group $\H^{i}(\spec K,F)$.  The (separable)
  \emph{index\/} of $\alpha$, denoted $\ind(\alpha)$, is $\gcd\{\deg
  L/K:\alpha|_{L}=0\}$ with $L$ a separable extension.
\end{defn}

In certain instances, the index of $\alpha$ is actually the
minimal degree of a field extension killing $\alpha$ (for example,
if $i=2$ and $F=\G_{m}$).

\begin{example}\label{E:brauer} When $F=\G_m$ and $i=2$, we can
  translate \ref{D:per-ind defn} into the language of non-commutative
  algebra.  A class $\alpha\in\H^2(\spec K,\G_m)=\Br(K)$ is
  represented by a central division $K$-algebra $D$.  The period of
  $\alpha$ is the order of $[D]$ in $\Br(K)$ and the index of $\alpha$
  is the square root of the rank of $D$.  This follows from basic
  theorems in non-commutative algebra \cite{farb-dennis} which
  describe the structure of central division algebras.
\end{example}

\begin{lem}\label{L:basic relation} Given $F$, $i$, and $\alpha$ as
  above,
  $\per(\alpha)|\ind(\alpha)$ and both have the same prime factors.
\end{lem}
\begin{proof} This is well-known, but we review the proof for the sake
  of completeness.  That $\per(\alpha)|\ind(\alpha)$ is a consequence
  of the existence, for any finite separable extension $f:\spec
  L\to\spec K$, of a trace map $f_{\ast}f^{\ast}F\to F$ (the
  ``corestriction'') such that the composition with the restriction is
  multiplication by $\deg L/K$.  Applying this to $L/K$ such that
  $\alpha_{L}=0$ shows that the degree of any such extension kills
  $\alpha$, whence the gcd kills $\alpha$.  To show that the prime
  factors are the same, suppose $L/K$ kills $\alpha$ and let $p$ be a
  prime number not dividing $\per(\alpha)$.  Suppose $L$ is Galois
  with group $G$.  Let $G_{p}$ be a $p$-Sylow subgroup of $G$ and
  $L_{p}$ the fixed field.  Thus, the degree of $L_{p}$ over $K$ has
  no factors of $p$.  We claim that $L_{p}$ kills $\alpha$.  Indeed,
  $\alpha_{L_{p}}$ is still $\per(\alpha)$-torsion and $p$ still acts
  invertibly on the subgroup generated by $\alpha$.  Thus, restricting
  further to $L$ and corestricting back to $L_{p}$, we see that
  $\alpha$ is trivial over $L_{p}$ if and only if it is trivial over
  $L$.  This removes $p$ from $\ind(\alpha)$.  The case of general
  separable $L/K$ follows from the case of Galois $L/K$.
\end{proof}

Thus, there is some minimal $\ell_{\alpha}$ so that
$\ind(\alpha)|\per(\alpha)^{\ell_{\alpha}}$.  The period-index problem
is to understand the behavior of the exponent $\ell_{\alpha}$ as $K$
and $\alpha$ vary.  In the context of elliptic curves, this was
considered by Lang-Tate \cite{lang-tate} and Lichtenbaum
\cite{lichtenbaum}.  Using \ref{E:brauer}, we see that the
period-index problem for the Brauer group is related to deciding how
large a division algebra is necessary to represent a given Brauer
class.  In this form, it was considered by Brauer in his foundational
paper and by Albert in subsequent work.

\begin{remark}
  Let $\pi:E\to \P^{1}$ be an elliptic surface over a finite field
  (with a section $\sigma$).  It is a classical result \cite[Exercise
  10.11]{silverman} that given an element
  $\alpha\in\Sha(k(\P^{1}),E_{\eta})$ of period $n$ prime to the
  characteristic of $k$, there is a field extension of degree $n$ of
  $k(\P^{1})$ killing $\alpha$.  In other words, for these particular
  classes of $\H^{1}(\spec k(\P^{1}),E_{\eta})$, the period equals the
  index.  Using Artin's isomorphism
  $\Br(E)\simto\Sha(k(\P^{1}),E_{\eta})$ (see \cite{grothbrauer3}), we
  see that this inadvertently also proves that
  $\per(\alpha)=\ind(\alpha)$ for any $\alpha\in\Br(E)$.
\end{remark}

\subsubsection{$\Br$ and $\Br'$}

The purpose of this paragraph is primarily to fix notation.  Given any
ringed topos $X$, there is a natural injection
$\Br(X)\inj\H^2(X,\G_m)$ (developed using the language of gerbes in
\S V.4 of \cite{giraud}).  The torsion subgroup of $\H^2(X,\G_m)$ is the
\emph{cohomological Brauer group\/}, denoted $\Br'$.

\begin{ques} When is the natural map $\Br\to\Br'$ an isomorphism?
\end{ques}

Grothendieck showed that $\Br=\Br'$ (where the latter is computed in
the \'etale topology) for Noetherian schemes of
dimension $1$ and regular schemes of dimension $2$
\cite{grothbrauer2}.  In his thesis \cite{gabber}, Gabber showed that
this is so for any scheme $X$ which can be written as a separated
union of two affines.  We give a proof of all of these results using
twisted sheaves in section \ref{S:gabber's theorems}.  A more recent
result of Gabber (mid-1990s, unpublished) established that $\Br=\Br'$
for any quasi-projective scheme, for which de Jong has given a twisted
proof \cite{dejong-gabber}.  There is an example of a non-separated
normal surface $Z$ such that $\Br(Z)[2]\neq\Br'(Z)[2]$
\cite{vistoli-kresch-etc.}.  The question is still quite open, even
for smooth non-quasi-projective threefolds over algebraically closed
fields!

\subsection{Geometry}

We follow the terminology and conventions of \cite{l-mb} when
discussing stacks and algebraic stacks.  Given a stack $\ms X$ over an
algebraic space $X$, we
will use $\ms I(\ms X)$ to denote the inertia stack, a representable
group-functor on $\ms X$ which associates to $\phi:T\to\ms X$ the
automorphisms of $\phi$.  (One can also explicitly describe the
inertia stack as $\ms I(\ms X)=\ms X\times_{\ms X\times\ms X}\ms X$, with
both maps being the diagonal morphism.  It is then an exercise to
check that this represents the sheaf of automorphisms, as described
above.  It appears that the inertia stack lacks a standard reference
in the literature.)  We will write $\sh(\ms X)$ to denote the
sheafification of $\ms X$, by which we mean the universal object among
sheaves $\ms F$ on $X$ admitting a map of stacks $\ms X\to \ms F$,
$\ms F$ being considered as a stack via the natural functor from sheaves
to stacks.

\subsubsection{Gerbes: generalities}

The canonical reference for gerbes and non-abelian cohomology is
\cite{giraud}, especially chapter IV.  Fix a positive integer $n$.
\begin{situation}\label{sit:prime to chars} 
  For the sake of simplicity, we will assume that $n$ is prime to the
  characteristics of $X$.
\end{situation}

\begin{defn} A \emph{gerbe\/} on $X$ is a stack $\ms X\to X$ such that
  the sheafification $\sh(\ms X)\to X$ is an isomorphism.  This is
  equivalent to the following two conditions:
  \begin{enumerate}
  \item[(i)] for every open set $U\to X$, there is a covering $V\to U$
    such that the fiber category $\ms X_{V}\neq\emptyset$;

  \item[(ii)] given an open set $U\to X$ and any two objects
    $x,y\in\ms X_{U}$, there is a covering $V\to U$ and an isomorphism
    $x_{V}\simto y_{V}$.
  \end{enumerate}
\end{defn}

\begin{defn} Given an abelian sheaf $A$ on $X$, an \emph{$A$-gerbe\/}
  on $X$ is a gerbe $\ms X$ along with an isomorphism $A_{\ms
    X}\simto\ms I(\ms X)$.  An \emph{isomorphism\/} of $A$-gerbes $\ms
  X$ and $\ms Y$ is a 1-morphism $f:\ms X\to\ms Y$ which is a
  1-isomorphism of stacks such that the natural map
$$A_{\ms X}\simto\ms I(\ms X)\to f^{\ast}\ms I(\ms Y)\simto A_{\ms X}$$
is the identity.
\end{defn}

The primordial example of an $A$-gerbe is the classifying stack
$\B{A}$ of $A$-torsors.  A general $A$-gerbe should then be thought of
as a bundle with fiber $\B{A}$ and ``$A$-linear'' gluing data.  The
original interest in these objects comes from the following result of
Giraud. For more on this subject, the reader is of course referred to
\S IV.3 of \cite{giraud}, and may consult \cite{mythesis} for a more
condensed discussion.

\begin{prop}[(Giraud)] There is a natural bijection between the set of
  isomorphism classes of $A$-gerbes and $\H^{2}(X,A)$.
\end{prop}

\begin{notn}
  We will write $[\ms X]$ for the class of $\ms X$ in $\H^2(X,A)$.
\end{notn}

Recall that given a sheaf $\ms F$ on a stack $\ms X$, there is a
natural right action $\ms F\times\ms I(\ms X)\to\ms F$.  Indeed, given
an object $U$ in the site of $\ms X$ and a section $(f,\sigma)$ of
$\ms F\times\ms I(\ms X)$ over $U$, we have an isomorphism
$\sigma^{\ast}:\ms F(U)\simto\ms F(U)$, and we let $(f,\sigma)$ map to
$\sigma^{\ast}f$.  On an $A$-gerbe, this becomes an action $\ms
F\times A_{\ms X}\to\ms F$.  

\begin{defn}\label{D:inertial action} Given an $A$-gerbe $\ms X$ and a
  sheaf on $\ms X$, the natural action $\ms F\times A\to\ms F$ will be
  called the \emph{inertial action\/}.
\end{defn}

\subsubsection{Gerbes: two special flavors}

We briefly describe two different types of $\m_{n}$-gerbes which will
be important for us.  Throughout, $X$ is assumed to be a $k$-scheme
for some field $k$ such that $n\in k^{\times}$.  (This last
requirement is not necessary if one is willing to work in the flat
topology.)

\begin{defn} A $\m_{n}$-gerbe $\ms X\to X$ is \emph{(geometrically)
    essentially trivial\/} if the class $[\ms X]$ has trivial image in
  $\H^{2}(X,\G_{m})$ (respectively, $\H^{2}(X\tensor_{k}\widebar
  k,\G_{m})$).
\end{defn}

\begin{defn}\label{D:gerbe of nth roots} Let $\ms M$ be an
  invertible sheaf on $X$.  The \emph{gerbe of $n$th roots\/} of $\ms
  M$, denoted $[\ms M]^{1/n}$, is the stack whose objects over $T$ are
  pairs $(\ms L,\phi)$, where $\ms L$ is an invertible sheaf on
  $X\times T$ and $\phi:\ms L^{\tensor n}\simto\ms M$ is an
  isomorphism.
\end{defn}
It is immediate that $[\ms M]^{1/n}$ is a $\m_{n}$-gerbe.  A
consideration of the Kummer sequence $1\to\m_n\to\G_m\to\G_m\to 1$
shows that every essentially trivial $\m_n$-gerbe is the gerbe of
$n$th roots of some invertible sheaf.  Geometrically essentially
trivial gerbes are more difficult to describe.

On the opposite side of the essentially trivial gerbes are the optimal
gerbes.

\begin{defn} A $\m_{n}$-gerbe $\ms X\to X$ is \emph{(geometrically)
    optimal\/} if the associated $\G_{m}$-gerbe has order $n$ in
  $\H^{2}(X,\G_{m})$ (resp.\ in $\H^{2}(X\tensor\widebar k,\G_{m})$).
\end{defn}
If $X$ is regular, this is the same as saying that the generic fiber
of $\ms X\to X$ has order $n$ in $\H^{2}(\eta_{X},\m_{n})$ (resp.\ in
$\H^2(\eta_{X\tensor\widebar k},\m_n)$), which is the same as saying
that $\ms X_{\eta}$ is the $\m_{n}$-gerbe associated to a division
algebra over $\kappa(\eta)$ of rank $n^{2}$.  If there is an Azumaya
algebra of degree $n$ on $X$ in the class $[\ms X]$, this is perhaps
saying that ``$\spec\ms A$ is an integral non-commutative surface
finite, flat, and unramified over its center.''

\section{Twisted sheaves and their moduli}\label{S:twisted sheaves and their moduli}

In this section, we will introduce the category of twisted sheaves and
show its suitability for studying certain types of questions about the
Brauer group.  We also give a summary of the results of
\cite{twisted-moduli}, which studies the properties of the moduli
spaces of twisted sheaves.

\subsection{Twisted sheaves}\label{S:twisted sheaves thing}

Let $X$ be an algebraic space.  We equip $X$ with a
reasonable flat topology: the big or small \'etale topology or the
fppf topology.  When we speak of a stack on $X$, we will mean a stack
in this topology.  When studying sheaves on a stack, we will use the
big \'etale or fppf topos.  As noted in the appendix, this requires
working with $\D(\qcoh)$ rather than the non-sensical
$\D_{\qcoh}(\Mod)$.  Throughout \S\ref{S:twisted sheaves thing}, $\ms
X\to X$ will be a fixed $D$-gerbe, where $D$ is a closed
subgroup of the multiplicative group $\G_m$.  We will write $C$ for
the dual of $D$; thus, $C$ is a quotient group of $\Z$.

\subsubsection{Twisted sheaves}

\begin{defn} An \emph{$\ms X$-twisted sheaf\/} is a sheaf $\ms F$ of
  (left) $\ms O_{\ms X}$-modules such that the natural inertial action
  $\ms F\times D_{\ms X}\to\ms F$ (see \ref{D:inertial action}) equals
  the right action associated to the left module action $D_{\ms
    X}\times\ms F\to\ms F$.
\end{defn}

Twisted vector bundles of rank $n$ naturally arise as ``universal
reductions of structure group'' for $\PGL_n$-bundles; the gerbe
classifying the twisting is precisely the coboundary of the
$\PGL_n$-bundle in $\H^2(\G_m)$ arising from the standard presentation
$1\to\G_m\to\GL_n\to\PGL_n\to 1$.  A similar $\m_n$-gerbe arises when
one studies reduction of structure group to $\SL_n$.  It is a helpful
exercise to write these gerbes down explicitly (as solutions to moduli
problems).  This may also be found in \cite{giraud}.

\begin{remark}\label{R:alternatives} 
  There are two other ways of describing twisted $\ms O$-modules which
  have arisen in the literature.  Both may be found in \Caldararu's
  thesis \cite{caldararu} (available electronically).  In case a
  faithful (i.e., nowhere zero) locally free twisted sheaf $\ms V$
  exists, one can apply the functor $\shom(\ms V,\ \cdot \ )$ and
  conclude (via a fibered form of Morita equivalence described in
  \cite{mythesis}) that the abelian category of twisted $\O$-modules
  is equivalent to the abelian category of right $\ms A$-modules,
  where $\ms A=\shom(\ms V,\ms V)$.

  Alternatively, one can give a cocyclic description of twisted
  sheaves (used by \Caldararu\ in his thesis \cite{caldararu},
  Yoshioka in his recent work \cite{yoshioka}, and others of a
  mathematico-physical bent) as follows: choose a hypercovering
  $U_{\bullet}\to X$ with a $2$-cocycle $a\in\m_{n}(U_{2})$
  representing the cohomology class $[\ms X]\in\H^{2}(X,\m_{n})$.
  (One can laboriously check \cite{mythesis} that indeed choosing such
  a setup is tantamount to choosing $\ms X$.)  A twisted sheaf is then
  a pair $(\ms F,\phi)$, with $\ms F$ a sheaf of $\O_{U_{0}}$-modules
  and $\phi:\pr_{1}^{\ast}\ms F\to\pr_{0}^{\ast}\ms F$ is an
  isomorphism on $U_{1}$ whose coboundary
  $\delta\phi=a\in\aut\pr_{00}^{\ast}\ms F$ on $U_{2}$.  In the cases
  where this definition has been applied to date, $X$ is usually
  quasi-projective over an affine, so it suffices to take a \v Cech
  hypercovering (by a theorem of Artin \cite{artin-joins}) in the
  \'etale topology.  (In fact, authors using this definition usually
  live over the complex numbers and use the classical topology.)
\end{remark}

\begin{para}
Since a $D$-gerbe is an Artin stack, there is a good theory of
quasi-coherent and coherent twisted sheaves.  Not surprisingly, this
theory is influenced by the representation theory of $D$.  
\end{para}

Given $\chi\in C$, for any quasi-coherent sheaf $\ms F$ on $\ms X$ there is a
$\chi$-eigensheaf $\ms F_{\chi}\subset\ms F$ corresponding to sections
where the action of the inertia stack $A_{\ms X}$ is via the character
$\chi:A\to\G_m$ and the right $\G_m$-action associated to the left
$\ms O$-module structure of $\ms F$.

\begin{prop}\label{P:q-coh splits} Suppose $\ms F$ is a
  quasi-coherent sheaf on $\ms X$.  The natural maps induce an
  isomorphism
$$\bigoplus_{\chi\in C(X)}\ms F_{\chi}\simto\ms F.$$
The eigensheaves $\ms F_{\chi}$ are quasi-coherent.
\end{prop}
\begin{proof}
  Using the specified identification $D_{\ms X}\cong\ms I(\ms X)$
  (which is part of the data describing a gerbe), the natural action
  of $\ms I(\ms X)$ on $\ms F$ makes $\ms F$ a sheaf of quasi-coherent
  comodules for the Hopf algebra structure on $\ms O_{D_{\ms X}}$.
  Thus, the proposition immediately follows from the basic
  representation theory of diagonalizable group schemes
  \cite{waterhouse}.
\end{proof}

Let $\pi:Y\to X$ be a quasi-compact $X$-space and $\ms
Y:=Y\times_{X}\ms X\to Y$ the pullback $D$-gerbe.  We will also denote
the morphism $\ms Y\to\ms X$ by $\pi$.

\begin{lem}\label{L:q-coh split is natural} If $\ms G$
  is a quasi-coherent sheaf on $\ms Y$, then the natural map
  $\pi_{\ast}\ms (\ms G_{\chi})\to\pi_{\ast}\ms G$ identifies
  $\pi_{\ast}(\ms G_{\chi})$ with $(\pi_{\ast}\ms G)_{\chi}$.
\end{lem}
\begin{proof}
As in the proof of \ref{P:q-coh splits}, we can view $\ms G$ as an 
$\ms O_{D_{\ms Y}}$-comodule.  In fact, the projection formula shows
that the comodule structure on $\ms G$ is the adjoint of the induced
$\ms O_{D_{\ms X}}$-comodule structure on $\pi_{\ast}\ms G$.  The
result follows.
\end{proof}

\begin{defn} An \emph{$m$-fold $\ms X$-twisted sheaf\/} is a sheaf
  $\ms F$ on $\ms X$ such that $\ms F=\ms F_{\chi_{m}}$, where
  $\chi_{m}:D\to\G_{m}$ is the $m$th power of the natural character
  $D\to\G_{m}$.
\end{defn}

Thus, $\ms X$-twisted sheaves are precisely the $1$-fold $\ms
X$-twisted sheaves.  Moreover, if $\ms X_n\to X$ represents the
cohomology class $n[\ms X]$, then one can check (e.g., using the
cocyclic formalism) that the category of $\ms X_n$-twisted sheaves is
(non-canonically) equivalent to the category of $n$-fold $\ms
X$-twisted sheaves.  The reason to introduce $m$-fold twisted sheaves
is for convenient bookkeeping, as the following lemma illustrates.

\begin{lem}\label{L:operations} 
  If $\ms F$ is $m$-fold twisted and $\ms G$ is $n$-fold twisted then
\begin{enumerate}
    \item all cohomology sheaves of $\ms F\ltensor\ms G$ are $(n+m)$-fold
    twisted;
    \item all cohomology sheaves of $\rshom(\ms F,\ms G)$ are
    $(n-m)$-fold twisted.
\end{enumerate}
If $\ms X$ is a $\m_{n}$-gerbe, then $n$-fold $\ms X$-twisted sheaves
are naturally identified with $0$-fold $\ms X$-twisted sheaves, which
are naturally identified via pushforward with sheaves on $X$.
\end{lem}
\begin{proof} This is practically a tautology.
\end{proof}

\begin{lem}\label{L:torsion criterion} Given a $\G_m$-gerbe $\ms Z\to
  X$, there is an invertible $n$-fold $\ms Z$-twisted sheaf if and
  only if $n[\ms Z]=0\in\H^{2}(X,\G_{m})$.  In particular, if there is
  a locally free $\ms Z$-twisted sheaf of rank $r>0$, then $r[\ms
  Z]=0\in\H^{2}(X,\G_{m})$.
\end{lem}
\begin{proof} The last statement follows by taking the determinant
  \ref{T:det} and applying the first statement.  To prove the first
  statement, we see from the fact that $\ms X_n$-twisted sheaves are
  identified with $n$-fold $\ms X$-twisted sheaves that it suffices to
  prove the statement assuming $n=1$.  It follows from standard
  techniques that the groupoid of invertible sheaves on $\ms X$ is
  equivalent to the groupoid of $1$-morphisms $\ms X\to\B{\G_m}$.
  (However, the reader should note that this is not the definition of
  $\B{\G_m}$, which only involves classifying invertible sheaves on
  algebraic spaces -- or perhaps only on affine schemes, depending
  upon one's ideological purity.)  The condition that an invertible
  sheaf $\ms L$ on $\ms X$ be $\ms X$-twisted is easily seen to
  correspond to the condition that the corresponding morphism
  $\phi_{\ms L}:\ms X\to\B{\G_{m,X}}$ induce the identity on
  inertia stacks via the specified morphisms $\ms I(\ms
  X)\simto\G_{m,\ms X}\to\phi_{\ms L}^{\ast}\ms I(\B{\G_m})$.  In other
  words, $\phi_{\ms L}$ is a morphism of $\G_m$-gerbes, hence an
  isomorphism, which means that the cohomology class of
  $\ms X$ is trivial.  Conversely, if $\ms X\cong\B{\G_{m,X}}$
  then there is an $\ms X$-twisted invertible sheaf arising from the
  canonical $1$-dimensional representation of $\G_m$.
\end{proof}

\begin{prop}\label{C:colimit} 
  If $X$ is Noetherian then a quasi-coherent $\ms X$-twisted sheaf is
  the colimit of its coherent $\ms X$-twisted subsheaves.
\end{prop}
\begin{proof} This follows from the fact that $\ms X$ is Noetherian,
  combined with standard results about quasi-coherent sheaves on
  Noetherian Artin stacks (\S 15 of \cite{l-mb}), along with the fact
  that any subsheaf of an $\ms X$-twisted sheaf is $\ms X$-twisted.
\end{proof}

\begin{para}
  We end this section with several results comparing various
  categories of twisted sheaves.
\end{para}
Giving a global section of $\pi:\ms X\to X$ yields a
``trivialization'' of the category of twisted sheaves.

\begin{lem}\label{L:section trivializes} 
  Given a section $\sigma:X\to\ms X$, the functor $\sigma^{\ast}$
  defines an equivalence of categories between the category of $\ms
  X$-twisted sheaves and the category of $\ms O_{X}$-modules.  This
  equivalence induces equivalences of the categories of quasi-coherent
  and coherent sheaves.
\end{lem}
\begin{proof} Given $\sigma$, there results an isomorphism
  $\widetilde\sigma:\ms X\simto\B{D}$ of $D$-gerbes on $X$
  which sends an object $\tau$ to $\isom(\tau,\sigma)$.  The natural
  character of $D$ gives an invertible $\B{D}$-twisted sheaf
  which pulls back to yield an invertible $\ms X$-twisted sheaf $\ms
  L$.  The functor $\ms W\mapsto(\pi^{\ast}\ms W)\tensor\ms L$ gives
  an equivalence between the category of $\ms O_X$-modules (with
  inverse given by $\ms V\mapsto\pi_{\ast}(\ms V\tensor\ms L^{\vee}$).

  The reader more comfortable with cocycles can also give a proof of this
  statement using cocyclic description of \ref{R:alternatives}.
\end{proof}

Using \ref{L:section trivializes}, we can prove that the category of
twisted sheaves is invariant under ``change of groups.''  Suppose
$\alpha\in\H^2(X,\m_n)$, and let $\iota_\ast\alpha$ be the image under
the natural map $\H^2(X,\m_n)\to\H^2(X,\G_m)$ induced by the inclusion
$\m_n\inj\G_m$.  There is a $\m_n$-gerbe $\ms Z\to X$ and a
$\G_m$-gerbe $\ms Y\to X$ with a $\m_n$-linear map
$\widetilde\iota:\ms Z\to\ms Y$ over $X$.

\begin{lem} The pullback functor associated to the morphism
  $\widetilde\iota$ yields an equivalence of the category of $\ms
  Y$-twisted sheaves with the category of $\ms Z$-twisted sheaves.
  This equivalence induces equivalences of categories of
  quasi-coherent and coherent sheaves.
\end{lem}
\begin{proof}
  Since $\ms Z$-twisted sheaves and $\ms Y$-twisted sheaves both form
  stacks on $X$, it suffices to prove the statement after covering $X$
  so that both gerbes become trivial.  Let $\ms L$ be an invertible
  $\ms Y$-twisted sheaf coming from a choice of section as in
  \ref{L:section trivializes}.  Given a sheaf $\ms W$ on $X$, there is
  a natural isomorphism $\widetilde\iota^{\ast}(\pi_{\ms Y}^{\ast}\ms
  W\tensor\ms L)\simto\pi_{\ms Z}^{\ast}\ms
  W\tensor\widetilde\iota^{\ast}\ms L$, and this immediately implies the result.

  One can also give a proof using the cocyclic description of twisted sheaves
  given in \ref{R:alternatives}: the image in $\G_m$ of the cocycle for $\m_n$
  corresponding to $\ms X$ gives the cocycle corresponding
  to $\ms Y$.  The cocyclic description of twisted sheaves now makes it
  clear that the categories are equivalent.
\end{proof}

\subsubsection{Twisted interpretations of $\Br=\Br'$ and period-index}

Using the technology of twisted sheaves, we can now clearly state both
the equality of $\Br$ and $\Br'$ and the period-index problem in
purely geometric terms.  We assume in this section that $X$ is
integral and Noetherian and that $D=\m_n$.

\begin{prop}\label{P:tw-interp}
  Let $\ms X\to X$ be a $\G_m$-gerbe and $m$ a positive integer.
  \begin{enumerate}
  \item $[\ms X]$ lies in $\Br(X)$ if and only if there
    is a non-zero locally free $\ms X$-twisted sheaf.
  \item If $X=\spec K$ is a reduced point, then $\ind [\ms X]$ divides
    $m$ if and only if there is a coherent $\ms X$-twisted sheaf of
    rank $m$.
  \item If $X$ is regular of dimension at most $2$, then $\ind [\ms
    X|_{\kappa(X)}]$ divides $m$ if and only if there is a locally free $\ms
    X$-twisted sheaf of rank $m$.
  \item The class $[\ms X]\in\H^2(X,\m_n)$ lies in the image of the
    coboundary map $\H^1(X,\PGL_n)\to\H^2(X,\m_n)$ if and only if
    there is a locally free $\ms X$-twisted sheaf of rank $n$ and
    trivial determinant (\ref{T:det}).
  \end{enumerate}
\end{prop}
\begin{proof}[Sketch of proof] The relevant fact is the following:
  there is a locally free $\ms X$-twisted sheaf of rank $r$ if and
  only if there is an Azumaya algebra of degree $r$ with class $[\ms
  X]$.  The equivalence arises by sending $\ms V$ to $\send(\ms V)$.
  The inverse comes by associating to an Azumaya algebra $\ms A$ its
  gerbe of trivializations: a section over a scheme $Y\to X$ is a
  locally free sheaf $\ms W$ along with an isomorphism of algebras
  $\send(\ms W)\simto\ms A$.  It is easy to check that this is a
  $\G_m$-gerbe with class $[\ms A]\in\Br(X)\subset\H^2(X,\G_m)$, and
  that the sheaves $\ms W$ patch together to yield a locally free
  twisted sheaf of rank $n$.  This is essentially carried out in \S
  V.4 of \cite{giraud}.  The third statement follows from the second
  one, combined with \ref{P:loc free on codimen 2}.

To prove the fourth statement, suppose $\ms V$ is a locally free $\ms
X$-twisted sheaf of rank $n$, so that $\send(\ms V)$ is the pullback
of an Azumaya algebra $\ms A$ of degree $n$ on $X$.  It is a simple
calculation to show that the class of $\ms A$ differs from $\ms X$ by
the image of $\det\ms V$ in $\H^2(X,\m_n)$ via the ``first Chern
class'' map.  Thus, the difference vanishes if and only if $\det\ms V$
is an $n$th power, in which case we can twisted $\ms V$ to produce a
locally free $\ms X$-twisted sheaf of rank $n$ with trivial determinant.
\end{proof}

Thus, these become problems about the existence of points in various
moduli spaces: locally free (or simply coherent) twisted sheaves (of a
given rank).  The period-index problem in particular will prove
amenable to this kind of analysis.  The question of $\Br$ and $\Br'$
is more profitably approached in terms of the $K$-theory of $\ms X$,
as we will describe in section \ref{S:gabber's theorems}.

\subsubsection{Elementary applications}\label{S:basic properties}

In this section, we assume that $X$ is Noetherian and that $D=\G_m$.

Here we investigate the applications of the theory we have developed
so far to the study of the Brauer group.  While seemingly vacuous, the
theory of twisted sheaves yields many of the basic results on the
Brauer group without requiring recourse to \'etale cohomology.  (The
skeptical reader should note that some of the fundamental properties
of gerbes and twisted sheaves do rely on the relationship to
cohomology, but this is where the similarity ends.  The point is that
we use sheaf-theoretic methods in place of explicit computations of
cohomological dimensions and the use of cohomological purity.)  The
essential result for all of these proofs, it turns out, is the
deceptively simple-looking \ref{C:colimit}.

\begin{lem}\label{L:extend my generic baby} Let $\eta\subset X$ be the
  scheme of generic points.  Any coherent $\ms X_{\eta}$-twisted sheaf
  extends to a coherent $\ms X$-twisted sheaf.
\end{lem}
\begin{proof} Let $\ms F$ be a coherent $\ms X_{\eta}$-twisted sheaf.
  The inclusion $\nu:\ms X_{\eta}\inj\ms X$ is quasi-compact, hence
  $\nu_{\ast}\ms F$ is a quasi-coherent $\ms X$-twisted sheaf.  The
  result follows from \ref{C:colimit}.
\end{proof}

\begin{lem}\label{L:non-zero coherent} There exists a non-zero
  coherent $\ms X$-twisted sheaf with support $\ms X$.
\end{lem}
\begin{proof} Over the reduced structure on the generic scheme of $X$,
  we have a $\G_{m}$-gerbe over the spectrum of a finite product of
  fields.  Thus, if we can produce a non-zero coherent twisted sheaf
  when $X$ is the spectrum of a field, we can push it forward to get
  such an object on the generic scheme of $X$ and then apply
  \ref{L:extend my generic baby}.  When $X$ is $\spec K$, any \'etale
  covering is finite over $X$.  Thus, there is a finite free morphism
  $Y\to X$ such that the gerbe $\ms Y:=Y\times_{X}\ms X$ has a
  section.  Once there is a section, we can apply \ref{L:section
    trivializes}.  Thus, there is a non-zero (in fact, locally free)
  twisted sheaf on $\ms Y$.  Pushing forward along the morphism $\ms
  Y\to\ms X$ yields a non-zero locally free $\ms X$-twisted sheaf by
  \ref{L:q-coh split is natural}. (For a generalization of this
  argument, see \ref{L:finite flat covering} below.)
\end{proof}

\begin{prop}\label{L:inj into gen pt} If $X$ is regular and integral
  with generic scheme
$\eta$,
then the restriction map $$\H^{2}(X,\G_{m})\to\H^{2}(\eta,\G_{m})$$ is
an injection.
\end{prop}
\begin{proof} Suppose $\ms X\times_{X}\eta$ represents the trivial
  cohomology class.  This means that there is an invertible $\ms
  X_{\eta}$-twisted sheaf $L_{\eta}$.  By \ref{C:colimit}, $L_{\eta}$
  has a coherent extension $L$ on all of $\ms X$.  In fact, there is a
  reflexive such extension.  On the other hand, this extension has
  rank 1 by construction.  But $X$, and therefore $\ms X$, is regular.
  As any reflexive module \emph{of rank 1\/} over a regular local ring
  (of arbitrary dimension) is free
  \cite[{\S}VII.4]{bourbaki-comm-alg}, we conclude that $L$ is an
  invertible twisted sheaf, whence $[\ms X]=0\in\H^{2}(X,\G_{m})$.
\end{proof}

\begin{cor} If $X$ is regular and quasi-compact then
  $\H^{2}(X,\G_{m})$ is torsion.
\end{cor}
\begin{proof} 
  This follows immediately from \ref{L:inj into gen pt} and the fact
  that the Galois cohomology of a field with abelian coefficients is torsion.
\end{proof}

Before proceeding, we recall a couple of lemmas which will also be
useful in section \ref{S:gabber's theorems}.

\begin{lem}\label{L:finite flat covering} If $Y\to X$ is a finite
  locally free covering and $\alpha\in\H^{2}(X,\G_{m})$, then there is
  a nowhere zero locally free twisted sheaf on $Y$ if and only if
  there is such a twisted sheaf on $X$.  If $X$ is quasi-compact, the
  same holds for locally free twisted sheaves of finite constant
  non-zero rank.
\end{lem}
\begin{proof} Fixing a gerbe $\ms X$ representing $\alpha$, we see
  that $\ms Y:=\ms X\times_{X}Y$ represents the pullback of $\alpha$
  to $Y$.  Furthermore, $\pi:\ms Y\to\ms X$ is a finite locally free
  morphism.  Thus, given any locally free twisted sheaf $F$ on $\ms
  Y$, $\pi_{\ast}F$ will be a locally free twisted sheaf on $\ms X$.
  It is clear that the rank of such a sheaf will be constant on a
  connected component.  If $X$ is quasi-compact then there are only
  finitely many connected components, whence one can arrange for the
  rank to be constant by taking appropriate direct sums on each component.
\end{proof}

\begin{lem}[(Gabber)]\label{L:local finite free cov} Given a local ring
  $A$ and a local-\'etale $A$-algebra $B$, there exists a finite free
  $A$-algebra $C$ such that for all maximal ideals $\mf m\subset C$,
  there is a map $B\to C_{\mf m}$ of $A$-algebras.
\end{lem}
\begin{proof}[Sketch of proof, following Gabber] By the local
  structure theory for \'etale morphisms \cite{raynaud-hensel}, there
  is a monic polynomial $f(x)=x^n+\sum_{i=0}^{n-1}a_ix^{n-i}\in A[x]$
  and a maximal ideal $\mf n\subset A[x]$ such that $f'\not\in\mf n$
  and $B\cong (A[x]/f(x))_{\mf n}$.  Let $C$ be
  $A[T_1,\ldots,T_n]/(\{\sigma_i(T_1,\ldots,T_n)-(-1)^ia_i;
  i=1,\ldots,n\})$, where $\sigma_i$ is the $i$th symmetric
  polynomial.  One can check that $C$ is a finite free $A$-module.  By
  construction, we have that $f(x)=\prod (x-T_i)$ in $C[x]$, so that
  there are factorizations $\phi_i:A[x]/(f(x))\to C$ over $A$ given by
  sending $x$ to $T_i$.  It is possible to see that $C$ is a finite
  free $A[x]/(f(x))$-module of rank $(n-1)!$ via $\phi_1$, from which
  it follows that there is some maximal ideal $\mf q\subset C$ lying
  over $\mf n\subset A[x]/(f(x))$, yielding a homomorphism $B\to
  C_{\mf q}$.  There is an obvious action of the symmetric group
  $\Sigma_n$ on $C$ over $A$, which is transitive on the closed fiber.
  Thus, composing the $A$-morphism $B\to C_{\mf q}$ with the
  automorphisms arising from the action of $\Sigma_n$, we see that
  there is in fact an $A$-map $B\to C_{\mf m}$ for every maximal ideal
  $\mf m\subset C$, as desired.
\end{proof}

\begin{prop}\label{L:loc free on dim 1} 
  If $X$ is semilocal of dimension at most 1 then there exists a
  locally free $\ms X$-twisted sheaf of positive rank.
\end{prop}
\begin{proof} When the cohomology class of $\ms X$ is assumed to be
  torsion, this will be subsumed by \ref{T:gabber hoobler twisted}
  below; thus, we only give a sketch here.

  We may assume that $X$ is connected.  It suffices to prove the
  statement assuming that $X$
  is reduced.  Indeed, the infinitesimal deformations of a locally free $\ms
  X_{\text{red}}$-twisted sheaf are easily seen to be unobstructed
  (using the fact that coherent cohomology vanishes on an affine
  scheme), so that a locally free $\ms X_{\text{red}}$-twisted sheaf
  may be deformed to a locally free $\ms X$-twisted sheaf of the same
  rank.

  If $X$ has dimension 0, then it is a disjoint union of reduced
  closed points.  Using the fact that the Brauer group and
  cohomological Brauer group coincide over fields (\S{X.5} of
  \cite{serre}) and elementary (e.g.\ Morita-theoretic)
  considerations, we see that there are locally free $\ms X$-twisted
  sheaves and that moreover there is at most one isomorphism class of
  $\ms X$-twisted sheaves with a given rank (at each point).
  (Briefly, the point is that by Morita theory one can
  see that two such differ by tensoring with an invertible sheaf,
  which is trivial over a semilocal ring.)

  Now suppose $X$ is local.  Applying \ref{L:finite flat covering} and
  \ref{L:local finite free cov}, we may replace $X$ by a finite
  semilocal extension $Y\to X$ such that $\ms Y:=\ms X\times_{X}Y$ is
  trivial over every localization of $Y$.  Choosing trivializations of
  the gerbe at each closed point of $Y$, it is immediate that there
  exists an integer $N$ such that for any closed point $y\in Y$, there
  is a locally free twisted sheaf on $\spec\O_{y,Y}$ of rank $N$.  By
  the previous paragraph, these are all isomorphic on the generic
  scheme of $Y$ (as $y$ varies).  It follows (since $Y$ is of
  dimension 1 at each closed point) that we can glue the local twisted
  sheaves to produce a locally free $\ms Y$-twisted sheaf.

  The general case proceeds as in the previous paragraph: the
  complement of the set of closed points of $X$ is open (since the
  dimension is at most $1$ at any point), so local models can be glued
  together to yield a locally free $\ms X$-twisted sheaf.
\end{proof}

\begin{cor}\label{C:torsiony} If $X$ is a semilocal 1-dimensional
  Noetherian scheme,
  $\H^{2}(X,\G_{m})$ is torsion and equals $\Br(X)$.
\end{cor}
\begin{proof} This is an application of \ref{P:tw-interp},
  \ref{L:torsion criterion}, and \ref{L:loc free on dim 1}.
\end{proof}

\begin{prop}\label{P:loc free on codimen 2} If $X$ is
  Noetherian there is a coherent twisted sheaf which is locally free
  at every point of codimension 1.  If $X$ is regular, there is a
  coherent twisted sheaf which is locally free at every point of
  codimension 2.
\end{prop}
\begin{proof} Using the first half of the proof of \ref{L:loc free on
    dim 1}, it is easy to find some open $V\subset X$ over which there
  is a locally free twisted sheaf $F$.  Suppose $p\in X\setminus V$
  has codimension 1 in $X$.  Let $i:\spec\O_{p,X}\inj X$ and $j:V\inj
  X$.  By \ref{L:loc free on dim 1}, we may assume (by taking a direct
  sum of $F$ with itself if necessary) that there exists a locally
  free twisted sheaf $G$ on $\spec\O_{p,X}$ such that
  $G_{\eta}=F_{\eta}$ (identifying the generic schemes of
  $\spec\O_{p,X}$ and $V$ with $\eta$ using $i$ and $j$,
  respectively).  Consider $Q:=i_{\ast}G\cap j_{\ast}F\subset
  F_{\eta}$.  This is a quasi-coherent twisted sheaf on $X$ which
  equals $F$ when pulled back by $j$ and $G$ when pulled back by $i$.
  Applying \ref{C:colimit} and the fact that localization commutes
  with colimits, we see that there is a coherent subsheaf $P\subset Q$
  such that $P|_{V}\cong F$ and $P|_{\spec\O_{p,X}}\cong G$. The locus
  where $P$ is locally free is thus an open set containing
  $V\cup\{p\}$.  By Noetherian induction, the first statement is
  proven.  (More concretely, there can only be finitely many
  codimension 1 points $p$ not in $V$, as they must be generic points
  of irreducible components of $X\setminus V$.) The second statement
  follows from the fact that any reflexive module over a regular local
  ring of dimension at most 2 is free.  Thus, the reflexive hull of
  any coherent twisted sheaf will be locally free in codimension 2.
\end{proof}

\begin{cor} If $X$ is regular of dimension at most 2 everywhere then
  the inclusion $\Br(X)\inj\H^{2}(X,\G_{m})$ is an isomorphism.
\end{cor}
\begin{proof} This follows from \ref{P:loc free on codimen 2} and
  \ref{L:torsion criterion}.
\end{proof}

\subsubsection{Gabber's theorems}\label{S:gabber's theorems}
We wish to show how twisted sheaves may be used to give an especially
streamlined proof of Gabber's theorem that $\Br=\Br'$ for affine
schemes.  (In fact, the proof given here also works for separated
unions of two affines, just as in Gabber's original result.)  Our
argument is a simplification of the argument of Hoobler \cite{hoobler}
which is itself a simplification of Gabber's proof.  By using twisted
sheaves rather than Azumaya algebras, the introduction of $K$-theory
into the proof becomes more transparent, as one can ``think in
modules'' from the beginning.  Our approach notably also avoids the
comparison of the ``Mayer-Vietoris sequence'' in non-abelian flat
cohomology with that in ordinary abelian flat cohomology, by absorbing
all of the cohomology (abelian and otherwise) into the underlying
gerbe $\ms X$.  The outline of our proof comes from Hoobler's paper
[\textit{ibid\/}.].

The main theorem of this section is the following.

\begin{thm}\label{T:gabber hoobler twisted} Let $X$ be an affine
  scheme and $\ms X$ an fppf $\m_{n}$-gerbe on $X$.  There exists an
  $\ms X$-twisted locally free sheaf of constant finite non-zero rank.
\end{thm}
\begin{cor}\label{C:gabber hoobler} If $X$ is an affine scheme then
  the natural injection $\Br(X)\inj\Br'(X)$ is an isomorphism.
\end{cor}
\begin{proof} Any torsion class $\alpha\in\H^{2}(X,\G_{m})$ (taken in
  the flat topology or \'etale topology) has a lift to a flat
  cohomology class in $\H^{2}(X,\m_{n})$.  The theorem gives a twisted
  vector bundle on a gerbe in any such class.  Taking its endomorphism
  algebra yields an Azumaya algebra with class $\alpha$.
\end{proof}

For the rest of this section, we assume that $X$ is an affine scheme
and $\ms X\to X$ is a $\m_n$-gerbe.

The proof of \ref{T:gabber hoobler twisted} follows from several basic
$K$-theoretic lemmas.  Let $R$ be a commutative unital ring and let
$K(R)$ denote the ring structure on the Grothendieck group of finite
projective $R$-modules of constant rank.  We recall some basic facts
from algebraic $K$-theory.

\begin{lem}\label{L:K-lemma} There is a natural map $\rho:K(R)\to\Z$
  determined by the rank function on modules.
  \begin{enumerate}
  \item The kernel of $\rho$ is a nil-ideal.
  \item Given projective modules $P$ and $Q$, if $[P]=[Q]$
    in $K(R)$ then there is $N>0$ such that $P^{\oplus
      N}\cong Q^{\oplus N}$.
  \item Given a class $\alpha\in K(R)$ of positive rank, some multiple
    $n\alpha$ has the form $[P]$ for a projective module $P$.
  \item If $P$ is a projective $R$-module such that $P^{\tensor n}$
    is free, then there is some $N>0$ such that $P^{\oplus N}$ is free.
  \end{enumerate}
\end{lem}
\begin{proof}[Sketch of proof.] Since the $K$-group classifies
  finitely generated projective modules, we immediately have that it
  commutes with colimits in $R$.  Thus, to prove any statement
  involving only finitely many elements it suffices to prove it
  assuming that $R$ is finitely generated over $\Z$.  In this case,
  $R$ has some Krull dimension $d$.

  A beautiful proof of the first statement is due to Gabber and may
  also be found in his thesis \cite{gabber}.  It is easy to see that
  $\ker\rho$ is generated by elements of the form $[P]-[R^{\oplus\rk
    P}]$, and thus that it suffices to prove the statement for such
  elements.  Now we can consider complexes $0\to R^{\oplus\rk P}\to
  P\to 0$, given by a morphism $\phi_i$, and it is easy to see that
  there are precisely $d$ such complexes such that for any point
  $p\in\spec R$, one of the complexes is exact at $p$ (i.e., one of
  the $\phi_i$ gives a trivialization of $P$ in a neighborhood of
  $p$).  But then the tensor product of the complexes
  $\phi_1\tensor\cdots\tensor\phi_d$ is exact, whence it vanishes in
  $K$-theory.  It follows that $([P]-[R^{\oplus\rk P}])^d=0$ in
  $K(R)$.

  The second statement is quite a bit more difficult, and relies on
  the Bass cancellation lemma.  (In the commutative case, it
  ultimately derives from a simpler statement due to Serre, which
  states that any projective module of rank at least $d+1$ has a
  non-zero free summand.)  This is beyond the scope of this paper; we
  refer the reader to \S 4 of Chapter IX of \cite{bass} for the
  details.  (The reader will also find a slightly different account of
  the rest of the $K$-theory we use here, including the other
  statements contained in this lemma.)

  The third statement may be proven as follows.  First, we may assume
  $R$ has finite Krull dimension $d$.  Scaling $\alpha$, we may assume
  that $\rk n\alpha\geq d$.  We can certainly write $n\alpha=[P]-[Q]$
  for some projective modules $P$ and $Q$.  Choosing a representation
  $R^N=Q\oplus Q'$, we find that $n\alpha+N=[P\oplus Q']$.  But
  $$\rk(P\oplus Q')=\rk(n\alpha+N)\geq d+N,$$ so we can find a free
  summand of $P\oplus Q'$ of rank $N$ by Serre's theorem.  The
  complementary summand yields a projective module representing
  $n\alpha$, as desired.

  To prove the final statement, let $r$ be the rank of $P$, and write
  $[P]=r+\eta\in K(R)$; we know from the first statement that $\eta$
  is nilpotent.  By hypothesis, $(r+\eta)^n\in\Z\subset K(R)$, and it
  follows from the nilpotence of $\eta$ that $n\eta=0$, so that
  $[P^{\oplus n}]=[R^{nr}]$.  Applying the second part, we conclude
  that there is $N$ such that $P^{\oplus N}$ is free, as required.
\end{proof}

\begin{cor}\label{C:k-theoretic hoobler tool} Given a faithfully
  projective $R$-module $P$ and a positive integer $n$, there exist
  non-zero free modules $F_{0}$ and $F_{1}$ and a faithfully
  projective $R$-module $\widebar P$ such that $P\tensor\widebar
  P^{\tensor n}\tensor F_{0}\cong F_{1}$.
\end{cor}
\begin{proof} In $K$-theory, the desired equality reads $a[P][\widebar
  P]^n=b$ for some positive integers $a$ and $b$.  In fact, applying
  \ref{L:K-lemma}(3), we see that it is enough to solve the equation
  for $[\widebar P]$ in $K(R)$.  Moreover, since we are allowed to
  scale by $a$ and $b$, it is clearly enough to solve the equation in
  $K(R)\tensor\Q$.  Writing $[P]=a+\beta$ with $\beta\in\ker\rho$, we
  see that since $\beta$ is nilpotent we can invert $(1/a)[P]$ and
  extract an $n$th root using formal power series with rational
  coefficients, thus yielding a solution in $K(R)\tensor\Q$.
\end{proof}

Combining \ref{L:local finite free cov} and \ref{L:finite flat
  covering}, we see that to prove \ref{T:gabber hoobler twisted} it is
enough to prove the following.  We will use the phrase ``locally
free'' to mean ``locally free of constant non-zero rank'' in what follows.

\begin{prop}\label{P:hoobler beef} Suppose that there is everywhere
  Zariski-locally on $X$ a locally free $\ms X$-twisted sheaf.  Then
  there is a global locally free $\ms X$-twisted sheaf.
\end{prop}
\begin{proof}  
  To simplify the notation, throughout this proof we will use the
  phrase ``twisted sheaf on $U$'' to mean ``$\ms X\times_X U$-twisted
  sheaf'' whenever $U\subset X$ is a subscheme.  It is enough to show
  that $$J=\{f\in A\ |\ \text{there is a locally free twisted sheaf on
    $\spec A_{f}$}\}$$ is an ideal.  The Zariski-local existence
  hypothesis shows that $J$ cannot be contained in any maximal ideal,
  hence if it is an ideal $J=A$.  (This is Quillen induction.) It is
  clear that $J$ is closed under multiplication by elements of $A$.
  To check that $J$ is closed under addition, we may reduce to the
  case where $X=U\cup V$ with $U$, $V$, and $U\cap V$ all affine, and
  there are locally free twisted sheaves on $U$ and $V$.  (Indeed,
  suppose $f,g\in J$; we wish to show that $f+g\in J$.  We may replace
  $X$ by $\spec A_{f+g}$ and thus assume that $f+g=1$, in which case
  $\spec A_f$ and $\spec A_g$ form an open cover of $X$.  Moreover,
  $\spec A_f\cap\spec A_g=\spec A_{fg}$ is affine.)

  Let $\ms P$ be a locally free twisted sheaf on $U$ and $\ms Q$ a
  locally free twisted sheaf on $V$; we may assume that $\rk\ms
  P=\rk\ms Q$.  Since $\ms X$ is a $\m_{n}$-gerbe, we see that $\ms
  P^{\tensor n}$ is naturally identified with a locally free untwisted
  sheaf (see \ref{L:operations}), and similarly for $\ms Q$.  By
  \ref{C:k-theoretic hoobler tool}, we see that there are non-zero
  locally free sheaves $\widebar P$ on $U$, $\widebar Q$ on $V$ and
  non-zero finite free modules $F_{0},F_{1}$ on $U$, $G_{0},G_{1}$ on
  $V$ such that
$$\ms P^{\tensor n}\tensor\widebar{P}^{\tensor n}\tensor F_{0}\cong
F_{1}$$ and $$\ms Q^{\tensor n}\tensor\widebar{Q}^{\tensor n}\tensor
G_{0}\cong G_{1}.$$ Thus, replacing $\ms P$ by $\ms
P\tensor\widebar{P}\tensor F_{0}$ and $\ms Q$ by $\ms Q\tensor\widebar
Q\tensor G_{0}$, we may assume that $\ms P^{\tensor n}$ and $\ms
Q^{\tensor n}$ are free modules on $U$ and $V$ respectively of the
same rank.  Now consider the situation on $U\cap V$.  Letting $P=\ms
Q^{\vee}\tensor\ms Q$ and $Q=\ms P\tensor\ms Q^{\vee}$, we have an
isomorphism of locally free twisted sheaves $\ms P\tensor P\cong\ms
Q\tensor Q$.  Furthermore, taking $n$th tensor powers, we see that
$P^{\tensor n}$ and $Q^{\tensor n}$ are both free (of the same rank).
Applying \ref{L:K-lemma}(4), we see that there is some $N>0$ such that
$P^{\oplus N}$ and $Q^{\oplus N}$ are free.  From the isomorphism $\ms
P\tensor P\cong\ms Q\tensor Q$ we deduce that $\ms P\tensor P^{\oplus
  N}\cong\ms Q\tensor Q^{\oplus N}$, from which it follows that $\ms
P^{M}\cong\ms Q^{M}$ for some integer $M>0$.  Thus, $\ms P^M$ and $\ms
Q^M$ are locally free twisted sheaves on $U$ and $V$ which glue on
$U\cap V$, as desired.
\end{proof}

Starting with the affine case \ref{C:gabber hoobler}, it is possible
to use twisted sheaves to prove the following.
\begin{prop}\label{C:gabber de jong} If $Y$ is a quasi-compact
  separated scheme admitting an ample invertible
  sheaf, then the natural injection $\Br(Y)\inj\Br'(Y)$ is an
  isomorphism.
\end{prop}
The reader is referred to \cite{dejong-gabber} for the proof.  The
idea is to start with a supply of twisted sheaves which are locally
free at selected points (using \ref{P:hoobler beef}), and to increase
the locus over which there exists a locally free twisted sheaf by
looking at kernels of general maps between the local models.  The
ample invertible sheaf enables one to make a Bertini argument when
studying a general such map between twisted sheaves (tensored with
powers of the ample invertible sheaf).

\subsection{Moduli: a summary}

In this section, we summarize the important aspects of the theory of
moduli of twisted sheaves, a subject which is treated in full detail
(with proofs) in \cite{twisted-moduli}.  In section
\ref{S:period-index}, we will use the structure theory of these moduli
spaces to prove results about the period-index problem for curves and
surfaces over finite, algebraically closed, and local fields.

\subsubsection{The moduli problem}

Let $X\to S$ be a smooth projective morphism (with a chosen relatively
ample $\O(1)$) with geometrically connected fibers, $n$ a positive integer
which is invertible on $S$, and $\ms X\to X$ a $\m_{n}$-gerbe.  Using
standard methods, one can prove the following (see
\cite{twisted-moduli}).

\begin{prop} The stack $\Tw_{\ms X/S}$ of $S$-flat torsion free $\ms
X$-twisted sheaves is an Artin stack locally of finite presentation over $S$.
\end{prop}
Given $T\to S$, the sections of $\Tw_{\ms X/S}$ over $T$ are $T$-flat
quasicoherent sheaves of finite presentation $\ms F$ on $\ms
X\times_{S}T$ which are $\ms X\times_{S}T$-twisted and such that for
every geometric point $t\to T$, the fiber $\ms F_{t}$ on $\ms X_{t}$
is torsion free.  (The reader is referred to \cite{twisted-moduli} for
basic results on associated points and torsion subsheaves on Artin
stacks.  Nothing surprising happens.)

As in the case of ordinary sheaves on a projective morphism, the stack
of all twisted sheaves is a nightmarish object.  For this and other
reasons, it is worthwhile to distinguish an open substack of
\emph{semistable\/} objects (and a further open substack of stable
objects).  In the cases at hand (namely for curves and surfaces), it
will turn out that this stack is in fact a GIT quotient stack;
however, the definition can be made in arbitrary dimension, and while
the stack has all of the properties one expects of a GIT quotient
stack, it is not clear if it is in fact such a quotient stack.

To define the stability condition, we will assume for the sake of
simplicity that $S=\spec k$ is the spectrum of a field.  In this case,
we can make use of the degree
function $\deg A_0(\ms X)\to\Q$, which can be defined as follows,
following \cite{vistoli}.  There is a rational Chow theory
$A_{\ast}(\ms X)$ for which there is a natural theory of Chern
classes.  Via proper pushforward, there is an isomorphism $A_0(\ms
X)\to A_0(X)\tensor\Q$.  

\begin{defn}
  The degree map $\deg:A_0(\ms X)\to\Q$ is the composition $A_0(\ms X)\to
  A_0(X)\tensor\Q\to\Q$, where the latter map is the degree in the
  usual Chow theory of $X$.
\end{defn}

In the following, we will write $\td_{\ms X}$ for the Todd class of
the tangent sheaf of $\ms X$.

\begin{defn}\label{D:polynomial} 
  Given a coherent sheaf $\ms F$ on $\ms X$, the \emph{geometric
    Hilbert polynomial of $\ms F$\/}, denoted $P_{\ms F}^{g}$, is the
  polynomial whose value at $m$ is $n\deg(\chern(\ms
  F(m))\cdot\td_{\ms X})$.
\end{defn}
It is a standard exercise to show that this in fact defines a (unique)
polynomial over $\Q$.  Furthermore, one can see that for torsion free
$\ms F$, the polynomial has degree equal to the dimension of $X$, with
leading coefficient equal to $\rk\ms F\deg_X\ms O(1)/(\dim X)!$.

\begin{remark}
  When the base $S$ is arbitrary, one can check that in fact the
  geometric Hilbert polynomial of an $S$-flat $\ms X$-twisted sheaf is
  constant in geometric fibers.  The reader is referred to 2.2.7.18 of
  \cite{twisted-moduli} for the details.
\end{remark}

With this additive function $K(\ms X)\to\Q[x]$ in hand, we define
(semi)stability in the usual way.  Recall that there is an ordering on
$\Q[x]$ given by the lexicographical ordering of coefficients.  In
this ordering, $f\leq g$ if and only if $f(a)\leq g(a)$ for all
sufficiently large (positive) integers $a$.

\begin{defn} With the above notation, a coherent sheaf $\ms F$ on $\ms
  X$ is \emph{semistable\/} if for all coherent proper subsheaves $\ms
  G\subset\ms F$, one has $(\rk\ms F)P^{g}_{\ms G}\leq(\rk\ms
  G)P^{g}_{\ms F}$.  The sheaf $\ms F$ is \emph{stable\/} if one has
  strict inequality for all coherent proper subsheaves $\ms G\subset\ms F$.
\end{defn}

\begin{remark}
  When the base $S$ is arbitrary, one can check that given an $S$-flat
  family $\ms F$ of coherent $\ms X$-twisted sheaves, there is an open
  subset of $S$ parametrizing semistable fibers $\ms F_s$ and a
  smaller open subset parametrizing geometrically stable fibers $\ms
  F_s$.  The reader is referred to 2.3.2.11 of \cite{twisted-moduli}.
\end{remark}

\begin{notn} The substack of semistable twisted sheaves is denoted
  $\Tw^{ss}\subset\Tw$; the substack of stable sheaves is
  $\Tw^{s}\subset\Tw^{ss}$.  The stack of semistable twisted sheaves
  with geometric Hilbert polynomial $P$ in the fibers is denoted
  $\Tw^{ss}_{\ms X/S}(P)$.
\end{notn}

The presence of the Brauer class of $\ms X$ can drastically simplify
the meaning of (semi)stability, as in the following.  We continue to
assume that $S=\spec k$ is the spectrum of a field.  Let $\Tw_{\ms
  X/S}(n)$ temporarily denote the (open) substack parametrizing
torsion free coherent $\ms X$-twisted sheaves of rank $n$.

\begin{lem}\label{L:stab-dumb}
  If $\ms X\to X$ is geometrically optimal, then the inclusion
  $\Tw_{\ms X/S}^{s}(n)\inj\Tw_{\ms X/S}(n)$ is an equality.
\end{lem}
\begin{proof}
  Given an $\ms X$-twisted sheaf $\ms V$ of rank $n$ on a geometric
  fiber of $X/S$, the fact that $\ms X\to X$ is geometrically optimal
  shows that any non-zero proper subsheaf $\ms W\subset\ms V$ must
  have full rank.  (I.e., $\ms V/\ms W$ is supported on a proper
  closed substack of $\ms X$.)  In this case, the inequality required
  for stability follows immediately from the fact that a non-zero
  coherent $\ms X$-twisted sheaf has non-vanishing geometric Hilbert
  polynomial.  (The reader is referred to 2.2.7.13 of
  \cite{twisted-moduli}; this seems to be a somewhat non-trivial fact,
  but it follows readily from the Riemann-Roch formula for
  representable morphisms of smooth Deligne-Mumford stacks with
  projective coarse spaces.)
\end{proof}

Just as with ordinary sheaves, one can fix the determinant of a
twisted sheaf.  As described in the appendix, there is a determinant
morphism $\Tw_{\ms X/S}\to\Pic_{\ms X/S}$ (more generally, if $X$ is
not smooth one can only define this morphism on the locus
parametrizing families whose fibers are perfect as objects in the
derived category of coherent twisted sheaves).  When the rank of $\ms
F$ is $n$, $\det\ms F$ is a $0$-fold twisted invertible sheaf (since
$\ms X$ is a $\m_n$-gerbe), and is thus identified via pushforward
with a section of $\Pic_{X/S}$.  This defines a map
$\Tw_{\ms X/S}(n)\to\Pic_{X/S}$.

\begin{defn}
  Given an invertible sheaf $\ms L$ on $X$ corresponding to
  $S\to\Pic_{X/S}$, the \emph{stack of twisted sheaves with
    determinant $\ms L$\/} is the stack-theoretic fiber product
  $\Tw_{\ms X/S}(n)\times_{\Pic_{X/S}}S$.  The stack of semistable $\ms
  X$-twisted sheaves of rank $n$, geometric Hilbert polynomial $P$,
  and determinant $\ms L$ will be denoted $\Tw^{ss}_{\ms X/S}(P,\ms
  L)$, and similarly for stable objects.
\end{defn}

\begin{prop} The stack $\Tw^{ss}_{\ms X/S}(P)$ has the following
  properties.
  \begin{enumerate}
  \item It is an Artin stack of finite presentation over $S$.
  \item It admits limits along any discrete valuation ring.
  \item The open substack $\Tw^{s}_{\ms X/S}(P)$ of stable twisted
    sheaves is a $\G_m$-gerbe over a separated algebraic space
    $\mTw^s_{\ms X/S}(P)$ of finite presentation over $S$.
  \item If $X$ is a relative curve or surface and $S$ is affine, then
    $\Tw^{ss}$ is a GIT quotient stack, and is thus corepresented in
    the category of schemes by a projective scheme.
  \end{enumerate}
  Given $\ms L\in\Pic_{X/S}(S)$, the stack $\Tw^{ss}_{\ms X/S}(P,\ms
  L)$ also satisfies i,ii, and iv, and the open substack $\Tw^s_{\ms
    X/S}(P,\ms L)$ is a $\m_n$-gerbe over a
  separated algebraic space of finite presentation, hence is itself separated.
\end{prop}

\begin{remark}
  As in the case of ordinary sheaves, the Brauer class corresponding
  to the $\G_{m}$-gerbe $\Tw^{s}\to\mTw^{s}$ is the famed Brauer
  obstruction to the existence of a tautological sheaf on $\mTw^{s}$.
\end{remark}

\begin{remark}
  One can set up a similar theory of semistable coherent sheaves for
  any proper smooth Deligne-Mumford stack $\ms X\to S$ with projective
  coarse moduli space $X$ along the lines indicated above.  From the
  point of view of period-index problems, this is particularly
  interesting when $\ms X$ is a $\m_n$-gerbe over an orbifold.  This
  point will be pursued elsewhere.
\end{remark}

With these abstract foundations in hand, we describe how the theory
specializes when $X$ is a curve or surface.

\subsubsection{Moduli of twisted sheaves on curves}\label{SS:mod-sh-c}

Fix a proper smooth curve $C$ over an algebraically closed field and a
$\m_n$-gerbe $\pi:\ms C\to C$.
It is a standard fact that $\H^{2}(C,\m_{n})=\Z/n\Z$ (with a canonical
generator being given by the first Chern class of the ideal sheaf of
any closed point of $C$).  Let
$\widebar{\delta}$ denote the (unique) fraction of the form $a/n$ such
that $0\leq a<n$ and $a\equiv[\ms C]\in\H^{2}(C,\m_{n})$.  

\begin{defn}
Given a sheaf $\ms F$ on $\ms C$, the \emph{degree\/} of $\ms
F$ is $n\deg c_1(\ms F)\in\Q$.  
\end{defn}

It is not hard to see that the denominator of the degree of $\ms F$ is
in fact a factor of $n$ and that when $\ms F$ is pulled back from $C$, the
degree computed above agrees with its degree as a sheaf on $C$.
(Without the correcting factor of $n$, the denominator would be a
factor of $n^2$, and degrees of pullbacks would be divided by $n$.)

\begin{notn} We let $\Tw_{\ms C/k}(r,d)$ denote the stack of twisted
  sheaves of rank $r$ and degree $d$ on $\ms C$.  The superscripts $s$
  and $ss$ will be added to denote the stable or semistable locus as
  above.  The coarse moduli space of $\Tw^s$ will be written as
  $\mTw^s$, as above.  When the determinant is identified with $\ms
  L$, we will write $\Tw_{\ms C/k}(r,\ms L)$ (with appropriate
  adornments for stable or semistable points).

  We let $\Sh_{C/k}(r,d)$ denote the stack of (untwisted) sheaves on
  $C$ with rank $r$ and degree $d$ and $\Sh_{C/k}(r,L)$ denote the stack of
  sheaves with fixed determinant $L$.  The (semi)stable locus is
  denoted with a superscript just as for twisted sheaves.
\end{notn}

\begin{prop}\label{P:on curves nada} 
  There is a non-canonical isomorphism $\Tw_{\ms
    C/k}(r,d)\simto\Sh_{C/k}(r,d-r\widebar\delta)$ which preserves the
  (semi)stable loci.
\end{prop}
\begin{proof}[Sketch of proof]
  The key to the proof is the observation that any $\ms C$-twisted
  invertible sheaf has degree congruent to $\widebar{\delta}$ modulo
  $\Z$.  By Tsen's theorem, there is an invertible $\ms C$-twisted
  sheaf, say $\ms T$, of degree $\widebar\delta$.  Just as in the
  proof of \ref{L:section trivializes}, the functor $\ms
  F\mapsto\pi_{\ast}(\ms F\tensor\ms T^{\vee})$ is an equivalence
  between the category of quasi-coherent $\ms C$-twisted sheaves and
  the category of quasi-coherent sheaves on $C$.  Moreover, given a locally
  free $\ms C$-twisted sheaf $\ms V$, it is easy to see that $\ms V$ is
  semistable (as a twisted sheaf) if and only if $\pi_{\ast}\ms
  (V\tensor\ms T^{\vee})$ is semistable (as a locally free sheaf on
  $C$).  Finally, if $\deg\ms V=d$, it follows from the usual
  formulas that $\deg\pi_{\ast}(\ms V\tensor\ms T^{\vee})=\deg(\ms
  V\tensor\ms T^{\vee})=d-r\widebar\delta$, as required.
\end{proof}

In particular, $\Tw^{ss}_{\ms C/k}(r,d)$ is a GIT stack (hence
corepresented by a projective variety).

The standard results concerning the moduli spaces of semistable
sheaves on smooth curves carry over to the twisted setting.  A
relatively exhaustive list of references for these classical results
may be found in \cite[Appendix 5C]{git}.  The final statement in the
following uses the main result of \cite{king-schofield}, which was not
available at the time \cite{git} was written.

\begin{cor}\label{P:curvy} Assume that $C$ has genus at
  least $2$. Fix an invertible $r$-fold $\ms C$-twisted sheaf $\ms L$
  of degree $d\in\Q$ (so that in particular $d-\widebar\delta
  r\in\Z$).
  \begin{enumerate}
  \item The stack $\Tw^{ss}_{\ms C/k}(r,\ms L)$ is smooth and is
    corepresented by a unirational projective variety 
    $\mTw^{ss}_{\ms C/k}(r,\ms L)$ of dimension $(r^{2}-1)(g-1)$.  The
    map on stable loci $\Tw^{s}_{\ms C/k}(r,\ms L)\to\mTw^s_{\ms
      C/k}(r,\ms L)$ is naturally a $\m_r$-gerbe.
  \item The stack $\Tw^{ss}_{\ms C/k}(r,d)$ is integral and smooth
    over $k$.  The stable locus is a 
    $\G_m$-gerbe over an algebraic space of dimension $r^2(g-1)+1$.
  \item If $d-r\widebar{\delta}$ and $r$ are relatively prime, then
    the open immersion $\Tw^{s}_{\ms C/k}(r,\ms L)\inj\Tw^{ss}_{\ms
      C/k}(r,\ms L)$ is an isomorphism.  In this case, $\mTw^{ss}_{\ms
      C/k}(r,\ms L)$ is a smooth rational projective variety
    representing the sheafification of the stack $\Tw^{ss}_{\ms
      C/k}(r,\ms L)$.  There is a tautological sheaf $\ms F$ on $\ms
    C\times\mTw^{ss}_{\ms C/k}(r,\ms L)$, and $\Pic(\mTw^{ss}_{\ms
      C/k}(r,\ms L))\cong\Z$.
  \end{enumerate}
\end{cor}

\begin{para} Now suppose the base field $k$ is not necessarily
  algebraically closed, and assume that $n$ is invertible in $k$.
  Write $S=\spec k$ and let $C\to S$ be a smooth projective
  geometrically connected curve with a section $p\in C(k)$.  There is
  a way to describe the spaces $\mTw^{ss}_{\ms C/k}(r,d)$ using the
  theory of Galois twists.  The Leray spectral sequence for $\G_{m}$
  and the choice of $p$ yield an isomorphism
  $\H^{2}(C,\G_{m})=\H^{2}(S,\G_{m})\oplus\H^{1}(S,\Pic_{C/k})$.
  Since $C$ has a point, it follows that
  $\H^{1}(S,\Pic_{C/k})=\H^{1}(S,\Pic^{0}_{C/k})$.  Similarly, there
  is a decomposition 
  \begin{equation}\label{Eq:only}
\H^{2}(C,\m_{n})=\H^{2}(S,\m_{n})\oplus
  \H^{1}(S,\Pic_{C/k}[n])\oplus\H^{0}(S,\R^{2}f_{\ast}\m_{n}).
\end{equation} 
The
  sheaf $\R^{2}f_{\ast}\m_{n}$ is in fact isomorphic to $\Z/n\Z$, and
  a splitting of the natural map
  $\H^{2}(C,\m_{n})\to\H^{0}(S,\R^{2}f_{\ast}\m_{n})=\Z/n\Z$ is given
  by sending $1\in\Z/n\Z$ to the gerbe $[\O(p)]^{1/n}$.

  In particular, the gerbe $\ms C$ gives rise to an element
  $\tau\in\H^{1}(S,\Pic_{C/k}[n])$ by projection.  Note that tensoring
  yields a map $\Pic_{C/k}[n]\inj\aut(\sh^{ss}_{C\tensor\widebar
    k/\widebar k}(n,\O(p)))$ (which is often an isomorphism
  \cite{pantev}).  By descent theory we see that varieties $V$ over
  $k$ which are geometrically isomorphic to $\mSh^{ss}_{C/k}(n,\O(p))$
  are classified (up to isomorphism) by
  $\H^{1}(S,\aut(\mSh^{ss}_{C/k}(n,\O(p))))$.  In particular, to any
  class $\tau\in\H^{1}(S,\Pic_{C/k}[n])$ is associated a twist of
  $\mSh^{ss}_{C/k}(n,\O(p))$.
\end{para}

\begin{prop}\label{P:twist identification} With notation as above, let
  $0\leq a<n$ correspond to the class of $\ms C\tensor\widebar k$ in
  $\Z/n\Z=\H^{0}(S,\R^{2}f_{\ast}\m_{n})$, and fix an invertible sheaf
  $\ms M\in\Pic(C)$.  Then $\mTw^{ss}_{\ms C/k}(n,\ms M)$ is the
  Galois twist of $\mSh^{ss}_{C/k}(n,\ms M(-ap))$ associated to the
  image of $[\ms C]$ in $\H^1(\spec k,\Pic_{C/k}[n])$ with respect to
  the decomposition of equation (\ref{Eq:only}) above.
\end{prop}
\begin{proof}[Sketch of proof]
  We first note that since $n$ is invertible in $k$ (by assumption),
  we see that the Brauer class of $\ms C$ will split over the
  separable closure of $k$, so that there is a $\ms C\tensor
  k^{\text{sep}}$-twisted invertible sheaf $\ms T$ as in \ref{P:on
    curves nada} such that $\ms T^{n}\cong\pi^{\ast}\ms O(ap)$ (so
  $\ms T$ has degree $a/n$).  The proof now 
  follows from the description of the isomorphism in \ref{P:on curves
    nada} combined with an explicit cocycle computation.
\end{proof}

In contrast to the classical case and \ref{P:curvy}(3), the existence of
a tautological family over $\mTw^{s}_{\ms C/k}(n,\ms M)$ (when $\deg\ms
M-r$ is prime to $n$) is by no means assured.  Since $\mTw^{s}$ is a
smooth projective rational variety with Picard group $\Z$, it follows
that the Brauer obstruction to the existence of a tautological sheaf
is in fact the pullback of a Brauer class from the base field.  As a
consequence, we see for example that the existence of a tautological
sheaf follows from the existence of a single stable twisted sheaf with
the given invariants over the base field.  (This in fact yields a
slightly different proof that the classical Brauer obstruction
vanishes in the untwisted case.)  One way to interpret this is as
follows: geometrically, one expects $\mTw^{s}$ to be quite similar to
the space of vector bundles.  The arithmetic of the twisting class is
responsible for the more complex structure of the projection
$\Tw^{s}\to\mTw^{s}$.  We will return to this issue later when we
study period-index problems for fibrations of curves.

\subsubsection{Moduli of twisted sheaves on surfaces}

Throughout this section, $\ms X\to X$ is a $\m_{n}$-gerbe on a smooth
projective surface over a field $k$ and $n$ is assumed prime to
$\ch(k)$.  The two moduli spaces which will be relevant are those
coming from essentially trivial gerbes and from optimal gerbes.  These
should be thought of as complementary cases: given an Azumaya algebra
$\ms A$ on $X$ with essentially trivial associated gerbe, one has that
the generic fiber $\ms A(\eta)$ is simply $\M_{n}(k(X))$.  If $\ms A$
has optimal associated gerbe, one has instead that $\ms A(\eta)$ is a
division algebra over the function field $k(X)$.  Twisted sheaves of
rank $n$ on an optimal $\m_{n}$-gerbe can be thought of as objects in
the ``right non-commutative Picard space'' of an integral
non-commutative surface finite over its center.  On the other hand,
moduli of essentially trivial twisted sheaves are best understood in
terms of a ``twisted stability condition'' on the moduli of ordinary
sheaves.  This is described in detail in \cite{twisted-moduli}.

The stack of twisted sheaves with fixed rank $n$, determinant $L$, and
geometric Hilbert polynomial $P$ will be denoted $\Tw(n,L,P)$.
Adornments, such as ${}^{ss}$, etc., will be added as necessary.  If
instead we fix the pair $(L,c_{2})$ of the determinant and degree of
the second Chern class, we will write $\Tw(n,L,c_{2})$.  (These two
ways of fixing the invariants are equivalent.)  The superscript $\mu$
will be used to denote the slope-stable locus, which is open in the
Gieseker-stable locus.

\begin{prop}\label{C:essentially trivial} Suppose $\ms X$ is
  geometrically essentially trivial.  There is a non-canonical
  isomorphism $$\Tw^{\mu}_{\ms X/k}(n,L,P)\tensor\widebar
  k\simto\Sh^{\mu}_{X/k}(n,L',Q)\tensor\widebar k,$$ with $Q$ an
  appropriate polynomial, where $\ms X=[L\tensor (L')^{\vee}]^{1/n}$.
\end{prop}

The proof is similar to the proof of \ref{P:on curves nada}.

\begin{remark}\label{R:c2-bd}
  In the above isomorphism, it is easy to see that the second Chern
  classes on the left and right sides differ by an amount which is
  bounded in terms of $L$, $\ms X$, and the intersection theory of
  $X$.  Indeed, if $\ms X\tensor\widebar k=[M]^{1/n}$, then the
  isomorphism arises by sending $\ms V$ to $\ms V\tensor\ms M^{\vee}$,
  where $\ms M$ is an $\ms X$-twisted invertible sheaf such that $\ms
  M^{\tensor n}\cong M$.  A standard computation shows that 
$$c_2(\ms V\tensor\ms M^{\vee})=c_2(\ms V)+(n-1)c_1(\ms
M^{\vee})c_1(\ms V)+n\binom{n}{2}c_1(\ms M^{\vee}).$$
Thus, if $\det\ms V$ is fixed, then $c_2(\ms V)$ and $c_2(\ms
V\tensor\ms M^{\vee})$ will differ by a fixed constant.
\end{remark}

In the case of an optimal gerbe, we have the following theorem.  The
interested reader will find a proof along with proof of several other
asymptotic properties (as $c_{2}$ grows) in \cite{twisted-moduli}.

\begin{thm}\label{T:irred} 
  Suppose $k$ is algebraically closed and $\ms X$ is optimal.  Given
  an invertible sheaf $L\in\Pic(X)$, there is a constant $A$ such that
  for all $c_{2}\geq A$, the stack $\Tw^{ss}_{\ms X/k}(n,L,c_{2})$ is
  integral, normal, and lci whenever it is non-empty.  Moreover, if
  $\Tw^{ss}_{\ms X/k}(n,L,c_2)$ is non-empty, then $\Tw^{ss}_{\ms
    X/k}(n,L,c_2+2r\ell)$ is non-empty for all $\ell\geq 0$.
\end{thm}

Using basic results from the appendix on elementary transforms, we
have the following non-emptiness result.

\begin{prop}\label{P:non-empty}
  Suppose $\ms X$ is geometrically optimal.  If
  there is a locally free $\ms X$-twisted sheaf of rank $n$ then for
  any integer $B$ and any invertible sheaf $L\in\Pic(X)$, the stack 
  $\Tw^{s}_{\ms X/k}(n,L,c_2)$ is non-empty for some $c_2>B$.
\end{prop}
\begin{proof}
  We may clearly assume $k$ is algebraically closed.  Since $\ms X$ is
  optimal, any torsion free $\ms X$-twisted sheaf of rank $n$ is
  automatically stable by \ref{L:stab-dumb}.  Moreover, if $\ms F$ is
  a torsion free $\ms X$-twisted sheaf of rank $n$ and $x\in X$ is a
  closed point, then the kernel $\ms F'$ of any surjection $\ms
  F\to\ms M_x$ onto a $\ms X\times_X x$-twisted sheaf of rank $1$ (a
  ``twisted skyscraper sheaf'') has the same determinant as $\ms F$
  and satisfies $c_2(\ms F')=c_2(\ms F)+2n$.  Thus, it suffices to
  show that there is a locally free $\ms X$-twisted sheaf of rank $n$
  and determinant $L$.  Suppose $\ms V$ is any locally free $\ms
  X$-twisted sheaf of rank $n$, and let $M=\det\ms V$.  If $\ms
  O_X(1)$ is a very ample invertible sheaf, then there is some $n$
  such that $L^{\vee}\tensor\det(\ms V(n))$ is very ample.  Let
  $D\subset X$ be a smooth divisor in the linear system
  $|L^{\vee}\tensor\det(\ms V(n))|$.  By Tsen's theorem, there is an
  invertible quotient of $\ms V(n)|_{\ms X\times_XD}$.  Applying
  \ref{C:tr-comp} yields the result.
\end{proof}

We note that recent results of Langer \cite{langer-castelnuovo} on moduli of
sheaves on arbitrary characteristic should give a method
for proving \ref{T:irred} in the absence of the optimality hypothesis.
In characteristic $0$, it is likely that existing methods of O'Grady
\cite{o'grady} will carry over without much change to the twisted
context, but we have not worked this out in detail.

\begin{remark}
  For the geometrically-minded, \ref{T:irred} plays an essential role
  in proving an algebraic analogue of results of Taubes on the stable
  topology of the space of self-dual connections on a fixed smooth
  $\PGL_{n}$-bundle over the 4-manifold underlying an algebraic
  surface.  We refer such readers to \cite{pgl-bundles}.
\end{remark}

The same methods used in the classical case also yield information
about the variation of the stack of stable twisted sheaves in a
family.  (The reader is referred to 3.2.4.25 of \cite{twisted-moduli}
for a proof of the following.)

\begin{prop}\label{P:dedekind family} Let $\ms X\to X\to S$ be a
  $\m_{n}$-gerbe on a smooth proper morphism over a locally Noetherian
  scheme with geometrically connected fibers of dimension 2, and
  assume that $n$ is invertible on $S$.  Suppose that for each
  geometric point $\widebar s\to S$, the fiber $\ms X_{\widebar x}\to
  X_{\widebar s}$ is optimal or essentially trivial.  Then the stack
  $\Tw^{s}_{X/S}(n,L,c_{2})\to S$ is a separated flat local complete
  intersection morphism for sufficiently large $c_{2}$.
\end{prop}

We point out that the hypothesis on the geometric fibers of $\ms X$ is
satisfied whenever $n$ is prime.  This usually suffices for the study
of period-index problems, as one can often reduce (by induction) to the case of
classes of prime period.

\section{Period-Index results}\label{S:period-index}

\subsection{Preliminaries}

\subsubsection{A cheap trick}\label{S:cheap trick} 
We mention here a simple trick which can be used to make certain base
extensions.  For the sake of concreteness, we only record the result
for field extensions.

\begin{prop}\label{P:cheap trick} Let $K$ be a field and
  $\alpha\in\Br(K)$ a class annihilated by $n$.  If $L/K$ is a finite
  field extension of degree $d$ and $n$ is relatively prime to $d$,
  then $\per(\alpha)=\per(\alpha|_L)$ and $\ind(\alpha)=\ind(\alpha|_L)$.
\end{prop}
\begin{proof} Write $X=\spec K$ and $Y=\spec L$.  The equality of the
  periods follows immediately from the fact that the map
  $\H^2(X,\m_n)\to\H^2(Y,\m_n)\to\H^2(X,\m_n)$ induced by pullback and
  trace is multiplication by $d$.  Suppose that $m$ is a positive
  integer with all prime factors dividing $n$.  If there is a locally free
  $\alpha_{Y}$-twisted sheaf on $Y$ of rank $m$ then pushing it
  forward to $X$ yields a locally free $\alpha$-twisted sheaf of rank
  $md$.  Taking endomorphisms yields a central simple $K$-algebra $A$
  of degree $md$ with $[A]=\alpha$.  If $D$ is the (unique up to
  isomorphism) central division $K$-algebra with class $\alpha$, we
  may write $A=\M_{r}(D)$ for some $r$.  Since $\per(\alpha)|n$, we
  see that $(\ind(\alpha),d)=1$, so the degree of $D$ is prime to $d$.
  Since $md=r\deg D$, we see that $d|r$.  Thus, $m=\ell\deg D$ for
  some $\ell$, and we conclude that $\M_{\ell}(D)$ is a central simple
  algebra of degree $m$ with class $\alpha$, and thus that
  $\ind(\alpha|_L)\geq\ind(\alpha)$.  Since the reverse inequality is
  obvious, the result follows.
\end{proof}

\begin{cor}\label{C:cheap corollary} 
  If $X$ is a $k$-scheme and $\alpha\in\H^{2}(X,\G_{m})$ has period
  $n$, then $\ind(\alpha)=n$ if and only if there is a field extension
  $k'\supset k$ of degree prime to $n$ such that
  $\ind(\alpha_{k'})|n$.
\end{cor}

\subsubsection{A reduction of geometric period-index problems to
  characteristic $0$}\label{S:reduction to char 0}

Let $k$ be an algebraically closed field, $K/k$ a finitely generated
field extension of finite transcendence degree and $\alpha\in\Br(K)$ a
Brauer class.  Let $W$ be the Witt ring of $k$; in
particular, $W$ is an absolutely unramified Henselian discrete
valuation ring with residue field $k$.

\begin{defn} Given a scheme $X$ and a point $x\in X$, an \emph{\'etale
    neighborhood of $x$\/} is an \'etale morphism $U\to X$ along with
  a lifting $x\to U$ over $x\to X$.
\end{defn}

\begin{prop}\label{P:lift} With the above notation, there exists an
  extension of discrete valuation rings $W\subset R$ such that
  \begin{enumerate}
  \item $R$ is essentially of finite type and formally smooth over $W$;

  \item there is an isomorphism $\rho$ between the residue field
    $\widebar R$ and $K$;

  \item for any finite extension $W\subset W'$, the base change
    $R\tensor W'$ is a discrete valuation ring with residue field
    isomorphic to $K$;

  \item there is a class $\widetilde{\alpha}\in\Br(R)$ with the same
    period as $\alpha$ and such that $\widetilde{\alpha}|_{\widebar
      R}$ corresponds to $\alpha$ via $\rho$.
  \end{enumerate}
\end{prop}
\begin{proof}
  The idea of the proof is to first fiber a proper model of $K$ by
  curves over an affine space, use the deformation theory of curves to
  lift the function field over $W$, and then use the deformation
  theory of Azumaya algebras to lift the Brauer class over an \'etale
  localization of the global model.

  First, we can choose a normal projective variety $X/k$ such that
  $k(X)\cong K$ (e.g., normalize the projective closure of an affine
  model!).  By Bertini's theorem, a generic hyperplane section of $X$
  will remain smooth in codimension $1$, so choosing a general pencil
  of very ample divisors on $X$, we can blow up $X$ to yield a
  birational model $\widetilde X\to\P^1$ with geometrically integral
  generic fiber of dimension $n-1$ which is smooth in codimension $1$.
  By induction, it follows that we can find a smooth curve $C$ over
  the function field $k(\P^{n-1})$ (with sections, if we like) such
  that the function field $K(C)$ is isomorphic to the original field
  $K$.

  Now consider the local ring $A$ of $\P^{n-1}_W$ at the generic point
  $\eta$ of the special fiber $\P^{n-1}_k$.  Over the residue field of
  $A$ we have a proper smooth curve $C$.  Standard deformation theory
  (as may be found in \cite{sga1}, for example) yields a proper smooth
  scheme $\mf C$ over the formal completion $\widehat A$ with
  geometrically connected fibers of dimension 1.  By Artin
  approximation (or Popescu's theorem), there results such a scheme
  over the Henselization of $A$.  The usual finite presentation
  arguments yield an \'etale neighborhood $U\to\P^{n-1}_W$ of $\eta$
  and a proper smooth relative curve $\ms C\to U$ such that the fiber
  over $\eta$ is the original curve $C$.  By shrinking $U$, we may
  assume that the fibers of $U$ over $\spec W$ are geometrically
  connected, hence that the fibers of $\ms C$ over $W$ are
  geometrically connected.  The same follows for any \'etale
  neighborhood of the generic point $\gamma$ of the special fiber of
  $\ms C$ over $W$.  Now $\ms O_{\ms C,\gamma}$ is a discrete
  valuation ring containing $W$ which satisfies all but possibly the
  last condition we require, and this remains true of any \'etale
  neighborhood of $\gamma$.  Thus, the Henselization of $\ms O_{\ms
    C,\gamma}$ is the colimit of subrings quasi-finite over $\ms
  O_{\ms C,\gamma}$ satisfying conditions 1 through 3.  But $\Br(\ms
  O_{\ms C,\gamma}^h)=\Br(K)$, so there is a lift of $\alpha$ to a
  class defined on $\ms O_{\ms C,\gamma}^h$, and this class descends
  to a class with the same period over some quasi-finite essentially
  finite type normal local extension $R$ of $\ms O_{\ms C,\gamma}$
  with residue field $K$.  This finishes the proof.
\end{proof}

Given a domain $B$, we temporarily write $F(B)$ for its field of fractions.

\begin{cor} With the notation of \ref{P:lift},
  $\ind\alpha=\ind\widetilde{\alpha}|_{R\tensor_W\widebar{F(W)}}$.
\end{cor}
\begin{proof} It is easy to see that there is a finite extension
  $W\subset W'$ such that
  $\ind\widetilde{\alpha}|_{R\tensor_WF(W')}=\ind\widetilde{\alpha}|_{R\tensor_W\widebar{F(W)}}$.
  But $R\tensor W'$ is a discrete valuation ring, so we can extend any
  division algebra over the fraction field $F(R\tensor
  W')=F(R)\tensor_{F(W)}F(W')$ to an Azumaya algebra over $R\tensor
  W'$ with Brauer class $\widetilde{\alpha}|_{W'}$.  Taking the
  residual algebra yields the result over the residue field of
  $R\tensor W'$, which is $K$.
\end{proof}

\subsection{Period and index on a fibration of curves}

\subsubsection{Relation between period-index and the existence of 
rational points in moduli}

\begin{lem}\label{L:good fibration} 
  Let $X$ be a smooth geometrically connected projective variety over
  a field $k$ of dimension $d>0$.  Suppose $X$ contains a $0$-cycle of
  degree $1$ (e.g., $k$ is finite).  Given any prime number $n$, there
  is a finite separable extension $k'\supset k$ of degree prime to $n$ and a
  birational equivalence of $X\tensor k'$ with a fibration $\widetilde
  X\to\P^{d-1}_{k'}$ with smooth generic fiber of dimension 1 and a
  rational point over $k'(\P^{d-1})$.
\end{lem}
\begin{proof} We may replace $k$ by the maximal separable prime to $n$
  extension and suppose that $k$ is infinite and $X$ has a rational
  point $p\in X(k)$.  Suppose first that $X$ is $\P^{d}$.  Taking the
  linear system of hyperplanes through $p$, one sees that the blowup
  $\bl_{p}\P^{d}$ fibers over $\P^{d-1}$ with a section and with
  generic fiber $\P^{1}_{k(\P^{d-1})}$.  Embedding $X$ in a projective
  space and choosing a generic linear projection yields a finite
  morphism $\nu:X\to\P^{d}$ which is \'etale at every point of
  $\nu^{-1}(\nu(p))$.  Now we can simply pull back the picture from
  $\P^{d}$.  The space $\widetilde X$ will then be the blowup of the
  fiber, which is just a finite set of reduced points.  Geometrically,
  this is the same as taking a general linear system of projective
  dimension $d-1$ containing $p$ in the base locus in any sufficiently
  ample complete linear system on $X$.
\end{proof}

\begin{lem} Let $X\leftarrow\widetilde X\to\P^{d-1}$ be as in
  \ref{L:good fibration}.  Given a class $\alpha\in\Br(X)$, one has
$$\ind(\alpha)=\ind(\alpha_{\widetilde X})=\ind(\alpha_{\widetilde
  X_{\eta}}),$$ where $\widetilde X_{\eta}$ is the generic fiber
curve.
\end{lem}
\begin{proof} This follows from the fact that $\widetilde X\to X$ is
  birational and the index is determined at the generic point.
\end{proof}
By construction, $\widetilde X_{\eta}$ has a rational point $p$.

Let $\pi:\widetilde X_{\eta}\to\eta=\spec K$ be a proper smooth
connected curve over a field with a rational point $p\in\widetilde
X_{\eta}(K)$.  Fix a class $\alpha\in\H^2(\widetilde X_{\eta},\m_n)$
and let $0\leq r<n$ be the class of $\alpha|_{\widebar K}$.  Let $\ms
C\to \widetilde X_{\eta}$ be a $\m_n$-gerbe representing $\alpha$.

\begin{prop}\label{P:rational points}  
  With the above notation, suppose the restriction of $\alpha$ to
  $p$ vanishes in $\Br(\kappa(p))$.  The following are
  equivalent:
\begin{enumerate}
    \item[(i)] $\per(\alpha)=\ind(\alpha)$;

    \item[(ii)] there is a locally free $\ms C$-twisted sheaf of
    rank $n$ and determinant of degree $r+1$;

\end{enumerate}
If $n$ is prime and the the $n$-cohomological dimension of $K$ is at
most $1$ then (i) and (ii) are equivalent to the existence of a
$K$-rational point on the coarse moduli space $\mTw^{s}_{\ms
  C/K}(n,\ms O((r+1)p))$, which is a smooth projective geometrically
rational variety with Picard number $1$.
\end{prop}
\begin{proof} The equivalence of (i) and (ii) in the absence of the
  degree requirement is just \ref{P:tw-interp}(3).  Thus, the content
  of the proposition lies in adjusting the determinant.  Suppose there
  is a locally free twisted sheaf $\ms F$ of rank $n$ on $\ms C$, and
  let $\ms M=\det\ms F$.  The assumption that
  $\alpha|_{p}=0\in\H^{2}(\spec\kappa(p),\G_{m})$ implies that there
  is a quotient of $\ms F|_{p}$ with any rank at most $n$.  Given such
  a quotient $Q$ of length $\ell$, it is easy to see that the
  determinant of the kernel $G$ of $\ms F\surj Q$ is $\det\ms F(-\ell
  p)$.  (This is just an elementary transform; see the appendix for
  details.)  By repeatedly applying this operation with various
  choices of $Q$ (and possibly dualizing to get the right sign), we
  can ensure that the determinant has degree $r+1$.  (In fact, we can
  achieve any chosen degree in this way.  The choice of $r+1$ will be
  made clear below.)  When $\operatorname{cd}_nK\leq 1$, we can use
  this argument at any effective Cartier divisor on $\widetilde
  X_{\eta}$ (as the $n$-torsion in the Brauer group of such a divisor
  vanishes) and therefore ensure that the determinant is precisely
  $\O((r+1)p)$.  The stability condition when $n$ is prime can be
  ensured in one of two ways: 1) if the image $\widetilde\alpha$ of
  $\alpha$ in $\Br(\widetilde X_{\eta})$ is $0$, then it is a
  classical fact that stable vector bundles of any given rank and
  determinant exist; 2) if $\widetilde\alpha\neq 0$, then $\alpha$ is
  optimal, so any locally free twisted sheaf of rank $n$ is stable and
  thus geometrically stable by numerical considerations.  (In general,
  only semistability is a geometric property, but in this case the two
  notions coincide.)  Finally, the cohomological dimension hypothesis
  ensures that any point of the coarse space $\mTw^s_{\ms C/K}$ will
  lift into the stack $\Tw^s_{\ms C/K}$, as the obstruction to doing
  so lies in $\Br(K)[n]$.

  The fact that $\mTw^s_{\ms C/K}(n,\ms O((r+1)p)$ is smooth,
  projective, and rational with Picard number $1$ follows from
  \ref{P:on curves nada}, the results cited above in section
  \ref{SS:mod-sh-c}, and the main results of \cite{king-schofield}.
\end{proof}

\subsubsection{An application: period and index on a geometric
  surface}
We give a proof of de Jong's theorem on the period-index problem for
surfaces over algebraically closed fields.

\begin{prop}\label{P:use de jong starr} Let $K$ be a field of
  transcendence degree 1 over an algebraically closed field and $C\to
  \spec K$ a smooth proper curve over $K$ with a rational point.
  Given a $\m_n$-gerbe $\ms C\to C$, let $0\leq r<n$ represent the
  class of $\ms C\tensor\widebar K$ in $\H^2(C\tensor\widebar
  K,\m_n)=\Z/n\Z$.  For any $\ms L\in\Pic(C)$ such that $\deg\ms L-r$
  is relatively prime to $n$, there exists a stable locally free $\ms
  C$-twisted sheaf of rank $n$ and determinant $\ms L$.
\end{prop}
\begin{proof} Since any Brauer class is defined over a finitely
  generated field over the prime field, we immediately reduce to the
  case where $K$ is finitely generated over its constant field.  By
  \ref{P:curvy}, $\mTw^s(n,\ms L)\tensor\widebar K$ is a non-empty
  smooth unirational projective variety.  (We cite the unirationality
  rather than rationality as it is far easier to prove.)  Applying the
  powerful theorem of Graber-Harris-Starr-de Jong (proven in
  \cite{starr-etc} in characteristic $0$ and generalized in
  \cite{dejong-starr} to positive characteristic), we conclude
  $\mTw^s(n,\ms L)$ has a rational point.  On the other hand,
  $\Tw^s(n,\ms L)\to\mTw^s(n,\ms L)$ is a $\m_n$-gerbe.  Since $K$ is
  $C_1$, any $\m_n$-gerbe over $K$ has a point, yielding a lift of the
  point in moduli to an object, as desired.
\end{proof}

\begin{lem}\label{L:ramification lemma} If $X$ is a smooth surface
  over an algebraically closed field $k$ and $\alpha\in\Br(k(X))$ has
  order prime to $\ch(k)$, then there is a birational morphism
  $\widetilde X\to X$ such that
  \begin{enumerate}
  \item the pullback of $\alpha$ to $\widetilde X$ has strict normal
    crossing ramification locus, and
  \item there is a fibration $\widetilde X\to\P^1$ with a section such that the
    restriction of $\alpha$ to the generic fiber is unramified.
  \end{enumerate}
\end{lem}
\begin{proof} It is an elementary fact that the ramification locus of
  $\alpha$ is pure of codimension $1$ (see \cite{artin-dejong},
  \cite{artin-mumford}, or \cite{saltman}).  Furthermore, if we blow
  up $X$, the only new ramification divisors which can appear must be
  the exceptional divisors.  The first statement follows by the
  embedded resolution of curves in surfaces.  To prove the second
  statement, we note that the ramification of $\alpha$ at a divisor
  $D$ is given by a cyclic extension of the function field $\kappa(D)$
  (\cite{artin-mumford},\cite{saltman}).  (In fact, a generically
  separable extension is sufficient for the present argument, and this
  follows easily from basic calculations in \'etale cohomology
  \cite{grothbrauer2} applied to the usual presentation
  for the sheaf of Cartier divisors.)  Let $E$ be the (normal
  crossing) ramification divisor of $\alpha$.  Choose a divisor $Z$ such that
  $Z$ intersects $E$ transversely at general points of irreducible
  components and $Z+E$ is very ample.  Choose a section $H$ of $|Z+E|$
  intersecting $Z$ transversely and intersecting $E$ transversely at
  general points of irreducible components over which the cyclic
  extensions measuring ramification are unramified.  It follows by an
  explicit cohomological calculation (or by the slightly more hands-on
  method of \cite{saltman}) that the blowup of $X$ at $H\cap{(Z+E)}$
  has no ramification on the exceptional divisors; i.e., the
  ramification of $\alpha$ on $\widetilde X$ is precisely contained in
  $E$.  This $\widetilde X$ is the total space of the pencil
  spanned by $H$ and $Z+E$, and we have just seen that $\alpha$ is
  unramified on the generic fiber.  The intersection points
  $H\cap(Z+E)$ give rise to sections of the pencil.
\end{proof}

\begin{thm}[(de Jong without restrictions)]\label{T:de jong per-ind} 
  If $K$ is a field of transcendence degree $2$ over an algebraically closed
  field $k$, then any $\alpha\in\Br(K)$ satisfies $\per(\alpha)=\ind(\alpha)$.
\end{thm}
\begin{proof} As above, we may assume that $K$ is
  finitely generated over $k$.  By the results of section
  \ref{S:reduction to char 0}, we may assume the base field has
  characteristic $0$.  (The reader can check that this proof in fact
  works whenever the period is prime to the characteristic exponent of
  the base field.)  Let $X$ be a smooth projective surface modeling
  the given function field and let $\alpha\in\Br(k(X))[n]$.  Applying
  \ref{L:ramification lemma}, we may replace $X$ with a fibration
  $\widetilde X\to\P^1$ such that $\alpha$ is unramified on the
  generic fiber.  Now we can apply \ref{P:use de jong starr}.  (The
  reader unfamiliar or uncomfortable with the ramification of a Brauer
  class can simply consider unramified classes -- these will suffice
  for the applications we give over finite fields -- at the expense of
  generality.)
\end{proof}

Elementary transforms can be used to refine \ref{T:de jong per-ind}.

\begin{cor}\label{C:det de jong} Let $X$ be a smooth projective
  surface over an algebraically closed field and $n$ an integer prime
  to $\ch(X)$.  The natural map $\H^1(X,\PGL_n)\to\H^2(X,\m_n)$ is
  surjective.
\end{cor}
\begin{proof} Given a $\m_n$-gerbe $\ms X\to X$, \ref{P:use de jong
    starr} shows that there is a locally free (stable) $\ms X$-twisted
  sheaf $\ms V$ of rank $n$.  Applying the methods used in the proof
  of \ref{P:non-empty} yields a $\ms V$ with trivial determinant.
  Taking $\send(\ms V)$ yields an Azumaya algebra whose cohomology
  class in $\H^2(X,\m_n)$ is precisely $[\ms X]$ (see \ref{P:tw-interp}(4))
\end{proof}

\begin{remark} One can also prove \ref{C:det de jong} for classes of
  order divisible by the characteristic (using the fppf topology to
  compute the cohomology) by applying the methods of \ref{P:non-empty}
  to fppf $\m_n$-gerbes on the surface $X$.
\end{remark}

\subsubsection{A local-to-global problem}\label{sec:local-global-problem}

A consequence of twisted methods is that one can solve the
period-index problem in general if one can characterize the points on
Galois twists of the moduli space $\sh^s_{C/k}(r,\ms L)$ of stable
vector bundles with fixed determinant of degree prime to the rank on a
smooth curve $C$ with a point.  Much is known about these varieties:
they are smooth and rational with Picard group $\Z$.  Moreover, one
knows that the canonical class of $\sh^s_{C/k}(r,\ms L)$ is twice an
anti-ample generator (so that $\sh^s$ is not a projective space!)
\cite{ramanan}.

It is tempting to ask the following question.

\begin{ques}\label{Conj:conj} Let $k$ be a field  
  and $C/k$ a smooth projective curve with a rational point.  Given an
  invertible sheaf $\ms L\in\Pic(C)$ and an integer $n$ prime to
  $\deg\ms L$, does any Galois twist of $\mSh^s_{C/k}(n,\ms L)$ coming from
  $\H^1(\spec k,\Pic^0_{C/k}[n])$ have a $k$-rational point?
\end{ques}

For a general curve $C/k$, the automorphism group of
$\mSh^s_{C/k}(n,\ms L)$ is precisely $\Pic^0_{C/k}[n]$, so the
conjecture is saying that in such cases, any form of
$\mSh^s_{C/k}(n,\ms L)$ should have a rational point.  This is closely
related to the period-index problem in the following way.

\begin{prop}\label{P:implication} 
  An affirmative answer to Question \ref{Conj:conj} implies that for
  any function field $K$ of dimension $d$ with
  algebraically closed constant field $k$ and any Brauer class
  $\alpha\in\Br(K)$, we have $\ind(\alpha)|\per(\alpha)^{2}$.
\end{prop}
\begin{proof}[Sketch of a proof] Let $n=\per(\alpha)$.  Given
  $K=k(X)$, we may birationally fiber $X$ over $\P^{d-1}$ such that
  \begin{enumerate}
  \item the generic fiber $C/k(x_1,\ldots,x_{d-1})$ is smooth with a
    section, and
  \item $\alpha$ is unramified on the generic fiber and trivial on the
    section.
  \end{enumerate}
  (Such a fibration comes from an appropriate linear system of
  hyperplane sections, and the section arises from a base point of the
  system.)  For an appropriate choice of $\ms L$, the coarse moduli
  space $\mTw^s(n,\ms L)$ has a rational point by \ref{P:twist
    identification} and \ref{Conj:conj}.  We have the following lemma.
\begin{lem}\label{L:univ-obs}
  Let $C/k$ be a proper smooth connected curve with a rational point
  $p$ and $\ms C\to C$ a $\m_n$-gerbe such that $\ms C\times_C p$ is
  trivial.  For any point $q\in\mTw^s_{\ms C/k}$, the index of the
  class in $\Br(\kappa(q))$ represented by $\Tw^s_{\ms
    C/k}\times_{\mTw^s_{\ms C/k}}q$ divides $n$.
\end{lem}
\begin{proof}
  Let $\ms F$ be the universal $\ms C$-twisted sheaf on $\ms
  C\times_k\Tw^s_{\ms C/k}$.  The universality of $\ms F$ implies that
  for any point $p\to\ms C$, the pullback $\ms F_p$ is $\Tw^s_{\ms
    C/k}$-twisted.  Thus, choosing $p:\spec k\to\ms C$ splitting the
  restriction $\ms C\times_C p$ yields a locally free $\Tw^s_{\ms
    C/k}$-twisted sheaf of rank $n$, as desired. 
\end{proof}
Thus, the obstruction to lifting the point into $\Tw^s$ has index dividing
$n$, so that, after making a base field extension of degree $n$, there
is a locally free twisted sheaf of rank $n$.  Pushing this sheaf
forward yields the result.
\end{proof}

On the other hand, Gabber's appendix to \cite{colliot-thelene}
(answering a question of Colliot-Th\'el\`ene in the affirmative) shows
that for a fixed prime $\ell$ and any positive integer $d$, there is a
smooth projective $d$-fold $X$ over $\C$ and an unramified Brauer
class $\alpha\in\Br(X)$ such that $\per(\alpha)=\ell$ and
$\ind(\alpha)=\ell^{d-1}$.  Taking a generic fibration of $X$ as in
\ref{P:implication}, we conclude the following.

\begin{cor}\label{C:0-cycle}
  Given a prime $\ell$ and a positive integer $d$, there is a smooth
  proper connected curve $C/\C(t_1,\ldots,t_{d-1})$ with a section and
  a $\m_{\ell}$-gerbe $\ms C\to C$ such that for every $L\in\Pic(C)$,
  the coarse moduli space $\mTw^s_{\ms C/k}(\ell,L)$ has the property
  that any closed point has degree divisible by $\ell^{d-3}$.
  Moreover, $\ms C\to C$ may be chosen so that its associated Brauer
  class is unramified on a proper smooth model of $C$ over $\C$.
\end{cor}
\begin{proof}
  We can fiber one of Gabber's examples to yield $\ms C\to C$
  representing a Brauer class of period $\ell$ and index $\ell^{d-1}$.
  (The ramification condition comes from the fact that the Brauer
  group is a birational invariant and Gabber's examples begin life as
  unramified classes.)  If there is a closed point whose degree has
  $\ell$-adic valuation smaller than $d-3$ then there is a locally
  free $\ms C$-twisted sheaf whose rank has $\ell$-adic valuation
  smaller than $d-1$.  (One additional factor of $\ell$ is potentially
  needed to lift the point into the moduli stack as in
  \ref{L:univ-obs}, and the other factor of $\ell$ comes from the fact
  that the stack parametrizes $\ms C$-twisted sheaves of rank $\ell$.)
  But this contradicts the fact that the index is $\ell^{d-1}$ (using
  \ref{C:cheap corollary}).
\end{proof}

There is a somewhat amusing consequence of this fact for
local-to-global problems over two-dimensional function fields.

\begin{prop}\label{P:loc-glob-prob}
For any positive integer $N$, there is a variety $V$ over
$\C(t_1,t_2)$ such that
\begin{enumerate}
\item $\dim V>N$;
\item $V$ is smooth, projective, geometrically connected, and
  geometrically rational;
\item $\Pic(V\tensor\widebar{\C(t_1,t_2)})\cong\Z$;
\item $V(\C(t_1,t_2)_\nu)\neq\emptyset$ for every non-trivial
  valuation $\nu$ on $\C(t_1,t_2)$, but
\item $V(\C(t_1,t_2))=\emptyset$.
\end{enumerate}
\end{prop}
\begin{proof}[Sketch of proof] 
  Choose $\ms C\to C$ as in \ref{C:0-cycle} with $d=3$, so that $C$ is
  the generic fiber of a proper fibration $\mc C\to\P^2_{\C}$ and $\ms
  C$ is the restriction of a $\m_n$-gerbe over $\mc C$.  Choose $L$ so
  that it is the restriction of an invertible sheaf on $\mc C$ and so
  that $\deg L-\ell\widebar\delta$ is relatively prime to $\ell$,
  where $\widebar\delta$ is as in \ref{SS:mod-sh-c}.  (After
  birationally modifying the base, we may assume that any given
  $L\in\Pic(C)$ extends, without doing any harm to the argument given
  here.)  Given a valuation $\nu$ of $\C(t_1,t_2)$, the valuation ring
  $\ms O_{\nu}$ will be centered at a point $p_{\nu}$ of $\P^2_{\C}$
  whose residue field has transcendence degree at most $1$.  The fiber
  $\mc C_{p_{\nu}}$ is a curve over the function field of a curve or
  over an algebraically closed field, so we may apply the methods of
  \ref{P:rational points} to find a locally free $\ms
  C_{p_{\nu}}$-twisted sheaf of rank $\ell$ with determinant
  $L_{p_{\nu}}$.  (When the fiber is not smooth, we can use the
  unobstructedness of locally free sheaves on curves to see that it
  suffices to solve the problem for the reduced structure on the
  fiber, then apply \ref{P:rational points} as written to the normalization
  and produce a twisted sheaf over the singular space by a gluing
  argument.)  It is a simple matter of deformation theory to see that
  this will deform over the complete local ring $\widehat{\ms
    O}_{\P^2,p_{\nu}}$.  Pulling back to the completion of
  $\C(t_1,t_2)$ with respect to $\nu$ yields a locally free $\ms
  C\tensor\C(t_1,t_2)_{\nu}$-twisted sheaf of rank $\ell$ and
  determinant $L$, which we may assume is stable (by a simple argument
  as in the proof of \ref{P:rational points}).  The reader can find
  further details of this kind of argument and applications of this
  method in more complicated (e.g., ramified) situations in
  \cite{per-ind-arith-surf}.

  On the other hand, a simple argument using the Leray spectral
  sequence shows that the map $\Br(\C(t_1,t_2))\to\Br(\mTw^s_{\ms
    C/\C(t_1,t_2)}(\ell,L))$ is an isomorphism.  It follows from the
  previous paragraph that the  
  universal obstruction $\Tw^s_{\ms
    C/\C(t_1,t_2)}(\ell,L)\to\mTw^s_{\ms C/\C(t_1,t_2)}(\ell,L)$ is
  split over every completion of $\C(t_1,t_2)$.  But any element of
  $\Br(\C(t_1,t_2))$ which is locally split everywhere must be split,
  and thus a $\C(t_1,t_2)$-rational point of $\mTw^s_{\ms C/\C(t_1,t_2)}(\ell,L)$
  would give rise to a $\ms C$-twisted sheaf of rank $\ell$,
  contradicting the fact that Brauer class of $\ms C$ has index $\ell^2$. 
\end{proof}

\begin{remark}
  The fact that the base field in \ref{P:loc-glob-prob} is
  $\C(t_1,t_2)$ is significant, because its unramified Brauer group is
  trivial.  It is quite a bit easier to make such examples over other
  function fields of transcendence degree $2$: if $P\to S$ is a
  Brauer-Severi scheme representing a non-trivial Brauer class over a
  proper smooth surface $S$, then it is easy to see that $P$ satisfies
  all of the conditions of \ref{P:loc-glob-prob}.
\end{remark}

\begin{remark}
  Using Th\'eor\`eme 10 of \cite{colliot-thelene}, one can make the
  counterexamples $V$ of \ref{C:0-cycle} be defined over $k(t_1,t_2)$ with
  $k$ any algebraically closed field.
\end{remark}

\subsection{Period and index via moduli on a surface}
Another way to reduce the period-index problem to a rationality
question is to use the known structure of the moduli spaces of twisted
sheaves on \emph{surfaces\/} along with classical estimates on the
existence of points (e.g., the Lang-Weil estimates for geometrically
integral varieties over finite fields).

\subsubsection{Period and index on a surface over a finite field}

\begin{thm}\label{T:per-ind finite surface} If $X$ is a proper
smooth geometrically connected surface over $\F_{q}$ and
$\alpha\in\Br(X)$ has period prime to $q$, then $\per(\alpha)=\ind(\alpha)$.
\end{thm}
For classes of period divisible by the characteristic, the methods
employed here cannot be applied (even though the moduli theory can be
developed).  Thus, (sadly) we can give no real insight into the
unramified problem for these classes.

\begin{remark}\label{R:p} If $K$ is the function field of a surface over a finite
  field of characteristic $p$, then it is easy to see that for
  $\alpha\in\Br(K)[p^{\infty}]$ we have $\ind(\alpha)|\per(\alpha)^2$.
  Indeed, the absolute Frobenius morphism $F:\spec K\to\spec K$ is
  finite free of degree $p^2$ and acts as multiplication by $p$ on
  $\Br(K)$.  Curiously enough, this result is \emph{stronger\/} than
  what we are able to prove in section \ref{S:ramified ramified} for
  ramified classes of period prime to the characteristic!
\end{remark}

\begin{remark}
  We thank Laurent Moret-Bailly for pointing out that the proof of
  \ref{T:per-ind finite surface} actually works for Brauer classes on
  a proper smooth geometrically connected surface over any perfect
  field $k$ with the property that any geometrically integral finite
  type $k$-scheme contains a $0$-cycle of degree $1$.  For example,
  the theorem applies to any PAC base field.  (Using \ref{P:strong
    pip} below, whose proof we have not included here, would allow one
  to remove the perfect hypothesis in the preceding sentence.) 
\end{remark}

Due to the inadequacy of the existing proofs of asymptotic properties
\cite{twisted-moduli} (e.g., applying only to geometrically optimal
classes and geometrically essentially trivial classes), we are forced
to make a slightly convoluted argument.  Let $PIP_{n}(k)$ denote the
phrase ``classes $\alpha$ of period dividing $n$ in the Brauer group
of any geometrically connected proper smooth surface over $k$ satisfy
$\per(\alpha)=\ind(\alpha)$.''

\begin{prop}\label{P:weak pip} 
  If $k$ is a perfect field and $PIP_{\ell}(L)$ for all finite
  extensions $L/k$ and all primes $\ell$ in a set of primes $S$ then
  $PIP_{n}(L)$ for all finite extensions $L$ and all $n$ in the
  submonoid $M$ of $\mathbf N$ generated by $S$.
\end{prop}
\begin{proof} Let $\alpha\in\Br(X)$ have period $n\in M$.  As above,
  we see that $n'\alpha$ has the period equal to index (where
  $n'=n/\ell$ for some prime factor $\ell$ of $n$).  Thus, $n'\alpha$
  has index $\ell$.  Since any division algebra over a field of index
  $d$ has a separable splitting field of degree $d$, we conclude that
  there is a generically \'etale finite map of normal surfaces $f:Y\to
  X$ of degree $\ell$ such that $f^{\ast}n'\alpha=0$.  By resolution
  of singularities in dimension $2$ (valid in any characteristic), we
  arrive at a regular proper surface $\widetilde Y$ and a morphism
  $\widetilde Y\to X$ which is generically finite of degree $\ell$.
  Since $k$ is perfect, we see that $\widetilde Y$ is smooth over $k$;
  in fact, setting $L=\H^{0}(Y,\ms O_{Y})$, we see that $\widetilde Y$
  is smooth, proper, and geometrically connected over $L$.  
  Applying $PIP_{n'}(L)$ to $Y$ and pushing forward at the generic
  point completes the proof.
\end{proof}

\begin{cor} To prove \ref{T:per-ind finite surface} it suffices to
  prove it for classes of prime order (prime to $q$).  Thus, it
  suffices to prove it for classes which are either geometrically
  optimal or geometrically essentially trivial.
\end{cor}

Using results of Artin and de Jong on general sections of an Azumaya
algebra, one can actually prove a stronger form of \ref{P:weak pip}.
Since it is not essential for our applications, we will only record
the statement.

\begin{prop}\label{P:strong pip} If $k$ is any field and
  $PIP_{\ell}(k)$ for all primes $\ell$ in a set of primes $S$ then
  $PIP_{n}$(k) for all $n$ in the submonoid of $\mathbf N$ generated by $S$.
\end{prop}  

\begin{proof}[Proof of \ref{T:per-ind finite surface} for
  geometrically optimal classes] Let $\ms X\to X$ be a geometrically
  optimal $\m_{n}$-gerbe.  Consider the stacks $\Tw^{s}_{\ms
    X/\F_q}(n,\ms O,b)$, which are defined over $\F_q$ since $\ms O$ is
  defined over $\F_q$.  By \ref{T:irred} and \ref{P:non-empty} (applied
  over $\widebar\F_q$, using \ref{T:de jong per-ind} as an input),
  there exists $b$ such that $\Tw^{s}_{\ms X/\F_q}(n,\ms O,b)$ is
  (non-empty and) geometrically integral. Thus, $\mTw^{s}_{\ms
    X/\F_q}(n,\ms O,b)$ is a geometrically irreducible (even
  projective) variety over $\F_{q}$, and furthermore since
  $\H^{2}(\spec \F_{q},\m_{n})=0$ we see that a rational point of
  $\mTw^{s}_{\ms X/\F_q}(n,\ms O,b)$ lifts to an object of
  $\Tw^{s}_{\ms X/\F_q}(n,\ms O,c_{2})$.  Thus, by the cheap trick
  \ref{S:cheap trick}, it suffices to find a rational point of
  $\mTw^{s}_{\ms X/\F_q}(n,\ms O,b)$ over $\F_{q}^{\text{non-$n$}}$,
  the maximal extension of degree prime to $n$.  But the Lang-Weil
  estimates \cite{lang-weil} show that any non-empty geometrically
  integral variety over $\F_{q}$ has rational points over
  $\F_{q}^{\text{non-$n$}}$.
\end{proof}

\begin{proof}[Proof of \ref{T:per-ind finite surface} for
  geometrically essentially trivial classes] In this case $\ms X\to X$
  is a $\m_{n}$-gerbe such that $\ms X\tensor\widebar{\F}_{q}$ admits
  an invertible twisted sheaf (has trivial Brauer class).  By
  \ref{C:essentially trivial}, the stack $\Tw^{\mu}(n,\ms O,c_{2})$ is
  geometrically isomorphic to $\Sh^{\mu}(n,\ms L,c_{2}')$ for some
  $\ms L$ (possibly defined only over an extension field, but this is
  unimportant) and a $c_{2}'$ which differs from $c_2$ by a fixed
  constant depending only upon $\ms X$ by \ref{R:c2-bd}.  On the other
  hand, the stack of semistable sheaves on $X$ is asymptotically
  geometrically irreducible and non-empty (in the sense of the
  statements of \ref{T:irred} and \ref{P:non-empty}), with $\Sh^{\mu}$
  as a dense open substack.  The proof in arbitrary characteristic
  (generalizing O'Grady's results in characteristic $0$) is contained
  in \cite{langer-castelnuovo}.  Applying the Lang-Weil estimates to
  the substack of $\mu$-stable points of $\Tw^{s}$ completes the
  proof.
\end{proof}

\subsubsection{A ramified result}\label{S:ramified ramified}

Using results of Saltman (\cite{saltman-fix} and \cite{saltman}), we
can extend our results to the ramified case.  We recall the relevant
results from Saltman's work (and we thank him for bringing this to our
attention).

\begin{prop}[(Saltman)]\label{L:saltman} Let $S$ be a regular scheme of
  dimension $2$ fibered over a field or discrete valuation ring
  containing a primitive $\ell$th root of unity, with $\ell$ a prime
  invertible on $S$.  Let $\alpha$ be a Brauer class at the generic
  point $\eta_{S}$ of prime order $\ell$.  There exist rational
  functions $f,g\in\kappa(S)$ such that $\alpha$ is unramified in
  $\kappa(S)(f^{1/\ell},g^{1/\ell})$.
\end{prop}
Saltman uses this to prove that $\ind|\per^2$ for Brauer classes in
the function field of a curve over a $p$-adic field.  Note that the
requirement that an $\ell$th root of unity exist is a minor one:
adjoining such a root creates a field extension of degree prime to
$\ell$, so using the methods of \ref{S:cheap trick} one concludes that
this has no effect on the period-index problem.  Similarly, the
requirement that $\alpha$ have prime order is immaterial when it comes
to period-index phenomena.

\begin{cor}\label{C:ramified finite} Let $X$ be a smooth projective
  geometrically connected surface over a finite field $k$. Any
  $\alpha\in\Br(\kappa(X))$ satisfies $\ind(\alpha)|\per(\alpha)^{3}$.
\end{cor}
\begin{proof} By \ref{R:p}, it suffices to prove this when the period
  of $\alpha$ is a prime $\ell$, and then we may assume $k$ contains
  all of the $\ell$th roots of unity by making a prime to $\ell$
  extension of the base field (which will not change the period or
  index by \ref{P:cheap trick}).  In this case, we can apply
  \ref{L:saltman} to find a finite extension of normal surfaces $Y\to
  X$ of degree $\ell^{2}$ such that $\alpha|_{Y}$ is unramified.
  Letting $L=\H^{0}(Y,\ms O_{Y})$, we can blow up $Y$ to arrive at a
  smooth geometrically connected surface $\widetilde Y$ over $L$
  equipped with a Brauer class $\alpha\in\Br(\widetilde Y)$ of order
  dividing $\ell$.  Since $Y$ is proper, $L$ is finite, and we may
  apply \ref{T:per-ind finite surface} to conclude that the index of
  $\alpha_{\widetilde Y}$ divides $\ell$.  Pushing forward a twisted
  sheaf of rank dividing $\ell$ along the morphism of generic points
  $\eta_{Y}\to\eta_{X}$ yields an $\alpha$-twisted sheaf of rank
  dividing $\ell^{3}$ on $\spec k(X)$ and taking the endomorphism ring
  yields a division algebra of degree dividing $\ell^{3}$ in the class
  $\alpha$, thus completing the proof.
\end{proof}

\subsubsection{Period and index on a surface over a local field}

Using \ref{T:per-ind finite surface} and \ref{P:dedekind family}, we
can prove a partial result on the period-index problem for surfaces
over local fields.  Throughout this section, $K$ denotes a local field
with integer ring $R$ and (finite) residue field $k$.

\begin{prop} Let $X$ be a proper smooth geometrically connected
  surface over $K$ which extends to a proper smooth relative surface $\mc
  X\to\spec R$.  If $\alpha\in\Br(X)$ has period prime to $\ch(k)$,
  then $\ind(\alpha)|\per(\alpha)^{2}$.  If $\alpha$ is unramified on
  $\mc X$, then $\ind(\alpha)=\per(\alpha)$.
\end{prop}
\begin{proof} First suppose $\alpha$ extends to all of $\mc X$.  On
  the special fiber, we have that $\per=\ind$.  It is not too
  difficult to show that one can construct an unobstructed locally
  free twisted sheaf of the desired rank using strict transforms (see
  e.g.\ 5.2.5 and 5.2.6 of \cite{h-l} for the untwisted version,
  \cite{mythesis} for the twisted version, or 3.2 of
  \cite{dejong-per-ind} for a treatment in terms of Azumaya algebras).
  Deforming it over the total space of $\mc X$ yields the result.

  If $\alpha$ is ramified along the special fiber, then extracting the
  $\per(\alpha)$th root of a uniformizer of $R$ will kill the
  ramification \cite[2.3.4]{artin-dejong} (which uses the main purity
  result of \cite{gabber}).  Thus, after making a finite free
  extension $R'\supset R$ of degree $\per(\alpha)$, we are reduced to
  the unramified case.  This is easily seen to imply the desired
  result.
\end{proof}

It would be nice to try to apply Saltman's techniques to the case of
local base fields $k$ and arrive at a result for ramified classes, but
it is not entirely clear how the property of smooth reduction behaves
under the construction of $\widetilde Y$ in the proof of
\ref{C:ramified finite}.

\begin{remark} It may be possible to extend these results in some form
  to the case of a surface with semistable or strictly semistable
  reduction.  To properly carry this out would require a study of
  semistable twisted sheaves and their deformations on semistable
  surfaces and their (possibly non-flat) infinitesimal thickenings.
  It is not particularly clear at the present time what exponent in
  the period-index relation one should expect in this case.
\end{remark}

\appendix

\section{A few facts about coherent sheaves on stacks}\label{S:appendix}
There is not much to recall here, as coherent sheaves behave just as
they do on schemes.  However, for the sake of completeness, we include
stacky versions of several well-known constructions and results.
First, it is worth making a brief remark about where precisely
quasi-coherent sheaves on a stack $\ms X$ live.  Of course, they are
quasi-coherent sheaves in one of the ringed topoi associated to $\ms
X$:

$$\ms X_{Zar}, \ms X_{fppf}, \ms X_{\textit{\'ET}}, 
\ms X_{\textit{lis-\'et}}, \ms X_{\textit{\'et}}.$$

The first three are all meant to be big topoi while the fourth is
developed in \cite{l-mb} -- but note that there are errors in their
treatment, which have been fixed by Olsson in \cite{olsson-sheaves}.
The fifth topos is the small \'etale topos: generated by \'etale
morphisms from schemes $U\to\ms X$.  Unless $\ms X$ is a
Deligne-Mumford stack, this is an unreasonable (mostly empty!)
choice.  We leave the proof of the following lemma to the reader (with
help from the references).

\begin{LEM} There are natural equivalences of abelian categories of
  quasi-coherent sheaves on $\ms X_{\text{Zar}}$, $\ms
  X_{\text{fppf}}$, $\ms X_{\text{\'ET}}$, and $\ms
  X_{\text{lis-\'et}}$.  When $\ms X$ is Deligne-Mumford, this
  equivalence extends to $\ms X_{\text{\'et}}$.
\end{LEM}

The derived category $\D(\qcoh(\ms X))$ has the same formal properties
when $\ms X$ is a stack as when it is a scheme.  In particular, if
$\ms X$ is regular of everywhere bounded homological dimension, any
coherent sheaf defines a perfect complex in $\D^{parf}(\qcoh(\ms
X))\subset \D(\qcoh(\ms X))$.  (Even when $\ms X$ is equidimensional
and quasi-compact, it is no longer true that the homological dimension
is related to the intrinsic dimension of $\ms X$, which can be
negative; it is however certainly bounded above by the dimension of a
smooth cover of $\ms X$.)

\begin{REMARK}
  Just as in the case of a scheme, the natural functor $\qcoh(\ms
  X)\to\Mod(\ms X)$ is not an exact functor of abelian categories when
  working with the big topologies, and therefore in the big \'etale
  and fppf topologies it is not reasonable to think about
  $\D_{\qcoh}(\ms X_{fppf})$ or $\D_{\qcoh}(\ms X_{\textit{\'ET}})$.
  It does make sense in the lisse-\'etale topology, but the entire
  setup requires great care in that context, as the formation of the
  lisse-\'etale topos is not functorial.  This is developed carefully
  and completely in \cite{olsson-sheaves}.  For Deligne-Mumford stacks
  one can work in the small \'etale topos $\ms X_{\textit{\'et}}$ and
  everything works just as for algebraic spaces.
\end{REMARK}

We recall the following theorem, essentially due to Knudsen and
Mumford \cite{mumford-knudsen} (to which the reader is referred for
further details).  Write $\operatorname{C}^{parf}(\ms X)$ for the
abelian category of perfect complexes of coherent sheaves and
$\operatorname{C}_{\cong}^{parf}(\ms X)$ for the subcategory
in which the morphisms are precisely the quasi-isomorphisms of perfect
complexes.  Write $\Pic_{\cong}(\ms X)$ for the Picard groupoid
(the category of invertible sheaves along with isomorphisms
between them on $\ms X$-schemes).

\begin{THM}\label{T:det} There is a natural additive and exact functor   
$$\det:\operatorname{C}^{parf}_{\cong}(\ms X)\to\Pic_{\cong}(\ms X)$$
which for every locally free sheaf $\ms V$ on $\ms X$, viewed as a
complex with $\ms V$ in degree $0$, takes the value
$\bigwedge^{\rk\ms V}\ms V$.
\end{THM}

The reader should refer to p.\ 23ff of [\emph{ibid\/}.] for the precise
properties the functor is assumed to have.  For our purposes, it
suffices to note that one can compute the determinant of a coherent
sheaf $\ms F$ on a regular stack by taking a locally free resolution
$\ms V^{\bullet}\to\ms F$ (when one exists) and tensoring the 
top wedge powers of the $\ms V^i$ with alternating signs.

\begin{DEFN} Given a coherent sheaf $\ms F$ on $\ms X$ whose image in
  $\D(\qcoh(\ms X))$ is perfect, the \emph{determinant of $\ms F$\/}
  is $\det\ms F$.
\end{DEFN}

\begin{PROP}\label{L:det} Let $\ms X$ be a regular Artin stack and
  $\ms C\subset \ms X$ an integral effective Cartier divisor.  A
  coherent sheaf $\ms F$ on $\ms C$ of generic rank $r$ has
  determinant $\ms O_{\ms X}(r\ms C)$ (as a perfect complex on $\ms X$).
\end{PROP}
\begin{proof} We would like to thank Moret-Bailly for pointing out a
  simplification of our original proof.  Let $\ms U=\ms X\setminus\ms
  C$ be the (open) complement of $\ms C$.  Since $\ms F|_{\ms U}=0$,
  there is an isomorphism $\sigma:\ms O_{\ms U}\simto\det\ms F|_{\ms
    U}$.  To show that $\sigma$ extends to an isomorphism $\ms O_{\ms
    X}(r\ms C)\simto\det\ms F$, it suffices to show this locally in
  the smooth topology of $\ms X$, as the descent datum on $\sigma$
  will descend the result back to $\ms X$.  Thus, it suffices to prove
  the result when $\ms X$ is a scheme.  To prove it in this case, note
  that an extension must exist for some $r$ by standard properties of
  invertible sheaves.  To determine $r$, it thus suffices to prove the
  result when $\ms X$ is the spectrum of a discrete valuation ring and
  $\ms C$ is the closed point, where the result is immediate.
\end{proof}

Let $\ms X$ be a regular Artin stack and $\iota:\ms C\inj\ms X$ an
integral effective Cartier divisor.  Suppose given a coherent sheaf
$\ms F$ on $\ms X$ and a coherent quotient of the restriction $\rho:\ms F|_{\ms
  C}\surj\ms Q$.
\begin{DEFN} With the above notation, the \emph{elementary transform
    of $\ms F$ along $\rho$\/} is
$$\ker(\ms F\to\iota_{\ast}\ms F_{\ms C}\to\iota_{\ast}\ms Q).$$
\end{DEFN}

\begin{COR}\label{C:tr-comp} 
  Given a coherent sheaf $\ms F$ on $\ms X$ and a quotient $\rho:\ms
  F_{\ms C}\surj\ms Q$, the determinant of the elementary transform of
  $\ms F$ along $\rho$ is $\det(\ms F)\tensor\O(-(\rk_{\ms C}\ms Q)\ms
  C)$.
\end{COR}
The corollary applies for example if $\ms X$ is quasi-compact and regular.

\def\cprime{$'$}


\begin{thebibliography}{10}

\bibitem{sga1}
{\em Rev\^etements \'etales et groupe fondamental}.
\newblock Springer-Verlag, Berlin, 1971.
\newblock S\'eminaire de G\'eom\'etrie Alg\'ebrique du Bois Marie 1960--1961
  (SGA 1), Dirig\'e par Alexandre Grothendieck. Augment\'e de deux expos\'es de
  M. Raynaud, Lecture Notes in Mathematics, Vol. 224.

\bibitem{caldararu}
Andrei {\Caldararu}.
\newblock Derived categories of twisted sheaves on {C}alabi-{Y}au manifolds,
  2000.
\newblock Thesis, Cornell University.

\bibitem{artin-joins}
M.~Artin.
\newblock On the joins of {H}ensel rings.
\newblock {\em Advances in Math.}, 7:282--296, 1971.

\bibitem{artinperind}
M.~Artin.
\newblock Local structure of maximal orders on surfaces.
\newblock In {\em Brauer groups in ring theory and algebraic geometry (Wilrijk,
  1981)}, volume 917 of {\em Lecture Notes in Math.}, pages 146--181. Springer,
  Berlin, 1982.

\bibitem{artin-dejong}
M.~Artin and A.~J. de~Jong.
\newblock Stable orders over surfaces, 2003.
\newblock Preprint.

\bibitem{artin-mumford}
M.~Artin and D.~Mumford.
\newblock Some elementary examples of unirational varieties which are not
  rational.
\newblock {\em Proc. London Math. Soc. (3)}, 25:75--95, 1972.

\bibitem{bass}
Hyman Bass.
\newblock {\em Algebraic {$K$}-theory}.
\newblock W. A. Benjamin, Inc., New York-Amsterdam, 1968.

\bibitem{bourbaki-comm-alg}
Nicolas Bourbaki.
\newblock {\em Commutative algebra. {C}hapters 1--7}.
\newblock Elements of Mathematics (Berlin). Springer-Verlag, Berlin, 1998.
\newblock Translated from the French, Reprint of the 1989 English translation.

\bibitem{castravet}
Ana-Maria Castravet.
\newblock Rational families of vector bundles on curves.
\newblock {\em Internat. J. Math.}, 15(1):13--45, 2004.

\bibitem{colliot-thelene}
Jean-Louis Colliot-Th{\'e}l{\`e}ne.
\newblock Exposant et indice d'alg\`ebres simples centrales non ramifi\'ees.
\newblock {\em Enseign. Math. (2)}, 48(1-2):127--146, 2002.
\newblock With an appendix by Ofer Gabber.

\bibitem{dejong-gabber}
A.~J. de~Jong.
\newblock A result of {G}abber, 2003.
\newblock Preprint.

\bibitem{dejong-per-ind}
A.~J. de~Jong.
\newblock The period-index problem for the {B}rauer group of an algebraic
  surface.
\newblock {\em Duke Math. J.}, 123(1):71--94, 2004.

\bibitem{dejong-starr}
A.~J. de~Jong and J.~Starr.
\newblock Every rationally connected variety over the function field of a curve
  has a rational point.
\newblock {\em Amer. J. Math.}, 125(3):567--580, 2003.

\bibitem{vistoli-kresch-etc.}
Dan Edidin, Brendan Hassett, Andrew Kresch, and Angelo Vistoli.
\newblock Brauer groups and quotient stacks.
\newblock {\em Amer. J. Math.}, 123(4):761--777, 2001.

\bibitem{farb-dennis}
Benson Farb and R.~Keith Dennis.
\newblock {\em Noncommutative algebra}, volume 144 of {\em Graduate Texts in
  Mathematics}.
\newblock Springer-Verlag, New York, 1993.

\bibitem{gabber}
Ofer Gabber.
\newblock Some theorems on {A}zumaya algebras.
\newblock In {\em The Brauer group (Sem., Les Plans-sur-Bex, 1980)}, volume 844
  of {\em Lecture Notes in Math.}, pages 129--209. Springer, Berlin, 1981.

\bibitem{giraud}
Jean Giraud.
\newblock {\em Cohomologie non ab\'elienne}.
\newblock Springer-Verlag, Berlin, 1971.
\newblock Die Grundlehren der mathematischen Wissenschaften, Band 179.

\bibitem{starr-etc}
Tom Graber, Joe Harris, and Jason Starr.
\newblock Families of rationally connected varieties.
\newblock {\em J. Amer. Math. Soc.}, 16(1):57--67 (electronic), 2003.

\bibitem{grothbrauer2}
Alexander Grothendieck.
\newblock Le groupe de {B}rauer. {II}. {T}h\'eorie cohomologique.
\newblock In {\em Dix Expos\'es sur la Cohomologie des Sch\'emas}, pages
  67--87. North-Holland, Amsterdam, 1968.

\bibitem{grothbrauer3}
Alexander Grothendieck.
\newblock Le groupe de {B}rauer. {III}. {E}xemples et compl\'ements.
\newblock In {\em Dix Expos\'es sur la Cohomologie des Sch\'emas}, pages
  88--188. North-Holland, Amsterdam, 1968.

\bibitem{guletskii}
V.~I. Guletski{\u\i} and V.~I. Yanchevski{\u\i}.
\newblock Torsion of the {B}rauer groups of curves defined over
  multidimensional local fields.
\newblock {\em Dokl. Nats. Akad. Nauk Belarusi}, 42(4):5--8, 121, 1998.

\bibitem{hoobler}
Raymond~T. Hoobler.
\newblock When is {${\rm Br}(X)={\rm Br}\sp{\prime} (X)$}?
\newblock In {\em Brauer groups in ring theory and algebraic geometry (Wilrijk,
  1981)}, volume 917 of {\em Lecture Notes in Math.}, pages 231--244. Springer,
  Berlin, 1982.

\bibitem{h-l}
Daniel Huybrechts and Manfred Lehn.
\newblock {\em The geometry of moduli spaces of sheaves}.
\newblock Aspects of Mathematics, E31. Friedr. Vieweg \& Sohn, Braunschweig,
  1997.

\bibitem{king-schofield}
Alastair King and Aidan Schofield.
\newblock Rationality of moduli of vector bundles on curves.
\newblock {\em Indag. Math. (N.S.)}, 10(4):519--535, 1999.

\bibitem{mumford-knudsen}
Finn~Faye Knudsen and David Mumford.
\newblock The projectivity of the moduli space of stable curves. {I}.
  {P}reliminaries on ``det'' and ``{D}iv''.
\newblock {\em Math. Scand.}, 39(1):19--55, 1976.

\bibitem{pantev}
Alexis Kouvidakis and Tony Pantev.
\newblock The automorphism group of the moduli space of semistable vector
  bundles.
\newblock {\em Math. Ann.}, 302(2):225--268, 1995.

\bibitem{kresch}
Andrew Kresch.
\newblock Hodge-theoretic obstruction to the existence of quaternion algebras.
\newblock {\em Bull. London Math. Soc.}, 35(1):109--116, 2003.

\bibitem{lang-tate}
Serge Lang and John Tate.
\newblock Principal homogeneous spaces over abelian varieties.
\newblock {\em Amer. J. Math.}, 80:659--684, 1958.

\bibitem{lang-weil}
Serge Lang and Andr\'e Weil.
\newblock Number of points of varieties in finite fields.
\newblock {\em Amer. J. Math.}, 76:819--827, 1954.

\bibitem{langer-castelnuovo}
Adrian Langer.
\newblock Moduli spaces and {C}astelnuovo-{M}umford regularity of sheaves on
  surfaces.
\newblock {\em Amer. J. Math.}, 128(2):373--417, 2006.

\bibitem{l-mb}
G{\'e}rard Laumon and Laurent Moret-Bailly.
\newblock {\em Champs alg\'ebriques}, volume~39 of {\em Ergebnisse der
  Mathematik und ihrer Grenzgebiete. 3. Folge. A Series of Modern Surveys in
  Mathematics [Results in Mathematics and Related Areas. 3rd Series. A Series
  of Modern Surveys in Mathematics]}.
\newblock Springer-Verlag, Berlin, 2000.

\bibitem{lichtenbaum}
Stephen Lichtenbaum.
\newblock The period-index problem for elliptic curves.
\newblock {\em Amer. J. Math.}, 90:1209--1223, 1968.

\bibitem{twisted-moduli}
Max Lieblich.
\newblock Moduli of twisted sheaves.
\newblock {\em Duke Math. J.}, 138(1):23--118, 2007.

\bibitem{mythesis}
Max Lieblich.
\newblock Moduli of twisted sheaves and generalized {A}zumaya algebras, 2004.
\newblock PhD Thesis, Massachusetts Institute of Technology.

\bibitem{pgl-bundles}
Max Lieblich.
\newblock Compactified moduli of projective bundles.
\newblock In preparation.

\bibitem{per-ind-arith-surf}
Max Lieblich.
\newblock Period and index in the Brauer group of an arithmetic
surface (with an appendix by Daniel Krashen), 2006.
\newblock Preprint.

\bibitem{milne}
James~S. Milne.
\newblock {\em \'{E}tale cohomology}, volume~33 of {\em Princeton Mathematical
  Series}.
\newblock Princeton University Press, Princeton, N.J., 1980.

\bibitem{git}
D.~Mumford, J.~Fogarty, and F.~Kirwan.
\newblock {\em Geometric invariant theory}, volume~34 of {\em Ergebnisse der
  Mathematik und ihrer Grenzgebiete (2) [Results in Mathematics and Related
  Areas (2)]}.
\newblock Springer-Verlag, Berlin, third edition, 1994.

\bibitem{nakayama}
T.~Nakayama.
\newblock \"Uber die direkte {Z}erlegung einer {D}ivisionsalgebra.
\newblock {\em Japanese Journal of Mathematics}, 12:65--70, 1935.

\bibitem{o'grady}
Kieran~G. O'Grady.
\newblock Moduli of vector bundles on projective surfaces: some basic results.
\newblock {\em Invent. Math.}, 123(1):141--207, 1996.

\bibitem{olsson-sheaves}
Martin Olsson.
\newblock Sheaves on {A}rtin stacks.
\newblock {\em J. Reine Angew. Math.}, to appear.

\bibitem{ramanan}
S.~Ramanan.
\newblock The moduli spaces of vector bundles over an algebraic curve.
\newblock {\em Math. Ann.}, 200:69--84, 1973.

\bibitem{raynaud-hensel}
Michel Raynaud.
\newblock {\em Anneaux locaux hens\'eliens}.
\newblock Lecture Notes in Mathematics, Vol. 169. Springer-Verlag, Berlin,
  1970.

\bibitem{saltman-fix}
David~J. Saltman.
\newblock Division algebras over {$p$}-adic curves.
\newblock {\em J. Ramanujan Math. Soc.}, 12(1):25--47, 1997.

\bibitem{saltman}
David~J. Saltman.
\newblock Correction to: ``{D}ivision algebras over {$p$}-adic curves''.
\newblock {\em J. Ramanujan Math. Soc.}, 13(2):125--129, 1998.

\bibitem{serre}
Jean-Pierre Serre.
\newblock {\em Local Fields}, volume 67 of {\em Graduate Texts in Mathematics}.
\newblock Springer-Verlag, New York, 1995
\newblock Second corrected printing of the 1979 original, translated
from the French by Marvin Jay Greenberg.

\bibitem{silverman}
Joseph~H. Silverman.
\newblock {\em The arithmetic of elliptic curves}, volume 106 of {\em Graduate
  Texts in Mathematics}.
\newblock Springer-Verlag, New York, 199?
\newblock Corrected reprint of the 1986 original.

\bibitem{vandenbergh}
M.~Van~den Bergh.
\newblock The algebraic index of a division algebra.
\newblock In {\em Ring theory (Antwerp, 1985)}, volume 1197 of {\em Lecture
  Notes in Math.}, pages 190--206. Springer, Berlin, 1986.

\bibitem{vistoli}
Angelo Vistoli.
\newblock Intersection theory on algebraic stacks and on their moduli spaces.
\newblock {\em Invent. Math.}, 97(3):613--670, 1989.

\bibitem{waterhouse}
William~C. Waterhouse.
\newblock {\em Introduction to affine group schemes}, volume~66 of {\em
  Graduate Texts in Mathematics}.
\newblock Springer-Verlag, New York, 1979.

\bibitem{yoshioka}
K{\=o}ta Yoshioka.
\newblock Twisted stability and {F}ourier-{M}ukai transform. {I}.
\newblock {\em Compositio Math.}, 138(3):261--288, 2003.

\end{thebibliography}
\end{document}